\documentclass[12pt]{article}
\usepackage{amsmath,amssymb,amsthm,amsfonts}
\usepackage{bbm}
\usepackage{graphicx} 
\usepackage[shortlabels]{enumitem}

\usepackage{subfigure,subcaption}
\usepackage{faktor}
\usepackage{mathtools}
\mathtoolsset{showonlyrefs}

\usepackage{hyperref}
\usepackage{svg}
\usepackage{enumitem}
\setlist{itemsep=3pt plus 1pt minus 1pt}

\usepackage{xcolor}

\renewcommand{\H}{\mathbb H}
\newcommand{\Z}{\mathbb Z}

\newcommand{\R}{\mathbb R}
\newcommand{\I}{\mathbb I}
\newcommand{\J}{\mathbb J}
\newcommand{\bs}{\backslash}
\newcommand{\FF}{\mathcal F}
\renewcommand{\S}{\mathcal S}

\newcommand{\minusone}{{-1}}

\newcommand{\pair}[2]{\left< #1, #2 \right>}
\newcommand{\geod}[2]{\overline{#1 #2}}

\newcommand{\numberthis}{\refstepcounter{equation}\tag{\theequation}}

\newtheorem{thm}{Theorem}[section]
\newtheorem{cor}[thm]{Corollary}
\newtheorem{prop}[thm]{Proposition}
\newtheorem{lem}[thm]{Lemma}

\theoremstyle{definition}
\newtheorem{defn}[thm]{Definition}

\theoremstyle{remark}
\newtheorem{rem}[thm]{Remark}

\newcommand{\dmo}{\DeclareMathOperator}
\dmo{\Id}{Id}
\dmo{\pt}{pt}
\dmo{\ab}{ab}
\dmo{\GL}{GL}
\dmo{\SO}{SO}
\dmo{\SU}{SU}
\dmo{\SL}{SL}
\dmo{\PSL}{PSL}
\dmo{\im}{Im}
\dmo{\re}{Re}
\dmo{\Ad}{Ad}
\dmo{\ad}{ad}
\dmo{\tr}{tr}
\dmo{\Log}{Log}
\dmo{\Homeo}{Homeo}
\dmo{\Vol}{Vol}
\dmo{\supp}{supp}
\dmo{\Isom}{Isom}

\newcommand{\N}{\mathcal{N}}

\newcommand{\norm}[1]{\left\| #1 \right\|}
\newcommand{\from}{\colon}
\newcommand{\sm}{\setminus}

\usepackage{stackengine,scalerel}
\stackMath

\newcommand\gtfrown{\mathrel{\raisebox{2pt}{\stackunder[1pt]{>}{\scaleto{\frown}{0.6ex}}}}}
\newcommand\ltfrown{\mathrel{\raisebox{2pt}{\stackunder[1pt]{<}{\scaleto{\frown}{0.6ex}}}}}

\makeatletter
\let\c@equation\c@thm
\makeatother
\numberwithin{equation}{section}

\title{Surface Subgroups for Cocompact Lattices of Isometries of $\H^{2n}$}
\author{Jeremy Kahn and Zhenghao Rao}
\date{}

\begin{document}
\maketitle
\begin{abstract}
    We prove the existence of surface subgroups within any cocompact lattice $\Gamma$ in $\mathrm{SO}(2n,1)$ for $n\geq2$. This result addresses the cases missing from the work of Hamenst\"adt in 2015, who constructed surface subgroups in cocompact lattices for all other rank-one simple Lie groups of non-compact type.
\end{abstract}


\tableofcontents
\section{Introduction}
A \emph{surface group},
for our purposes, 
is the fundamental group of a closed hyperbolic surface.
A surface \emph{subgroup} of a group is a subgroup that is a surface group.

Surface subgroups were found in \emph{arithmetic} cocompact lattices of $\Isom(\H^3)$ by Lackenby \cite{Lac10},
and in non-uniform lattices in $\Isom(\H^3)$ by Masters-Zhang \cite{MZ08,MZ09} and Baker-Cooper \cite{BC15}.

In \cite{KM12b},
Kahn and Markovic constructed ubiquitous nearly geodesic surface subgroups of all cocompact lattices in $\Isom(\H^3)$, 
by building a $\pi_1$-injective immersed surface (in the quotient manifold) out of immersed pairs of pants.
This result was used by Cooper and Futer \cite{CF19} to build ubiquitous surface subgroups in a given non-uniform lattice in $\Isom(\H^3)$; Kahn and Wright \cite{KW21} later proved the existence of nearly geodesic surface subgroups in the same setting. In \cite{Ham15}, Hamenst\"adt constructed surface subgroups of cocompact lattices in any rank-1 simple Lie group \emph{except} for $\Isom(\H^{2n})$ for $n \ge 2$. Kahn, Labourie and Mozes \cite{KLM24} built surface subgroups in cocompact lattices in center-free, complex semisimple Lie group of non-compact type.
In this paper, we show that for all $n\geq 2$, every cocompact lattice of isometries of $\H^{2n}$ has a surface subgroup. 
Here is our main result.

\begin{thm}\label{main thm}
    Suppose $\Gamma<\SO^+(2n,1)$ is a cocompact lattice with $n>1$. Then $\Gamma$ contains surface subgroups.
\end{thm}

Combining this with Hamenst\"adt's theorem yields a complete answer to the existence of surface subgroups in uniform lattices of simple rank-one Lie groups. Theorem \ref{main thm} also presents more positive concrete examples to Gromov's question on the existence of surface subgroups in one-ended hyperbolic groups. Although we do not state a quantitative result, the surface subgroups we construct here could be made $K$-quasi-Fuchsian for any $K>1$, as mentioned by Hamenst\"adt in \cite{Ham15}.

In this introduction,
we describe the basic idea used in \cite{KM12b} (and \cite{Ham15}, \cite{KW21}, \cite{KLM24}, \cite{Rao25}),
then explain the difficulty of applying it to $\H^{2n}$,
and finally outline our approach to resolve the issue.

In \cite{KM12b},
the authors consider ``good pants'' in $M = \Gamma \bs \H^3$,
which are $\epsilon$-nearly isometric immersions of the ``$R$-perfect pants'' into $M$;
the $R$-perfect pants is the unique hyperbolic pair of pants with three geodesic boundaries of length $R$. Each good pants in $M$ has three \emph{short orthogeodesics}, each connecting two of the boundary cuffs; each short orthogeodesic has two \emph{feet}, which are unit normal vectors to the boundary geodesics, pointing toward the orthogeodesic.
We then consider the distribution of the feet of good pants around each ``good'' closed geodesic $\gamma$ in $M$. We see that these are evenly distributed in the unit normal bundle $N^1(\gamma)$, and this implies that there is a permutation $\sigma\from \FF_\gamma \to \FF_\gamma$ (where $\FF_\gamma$ is the set of feet of good pants on $\gamma$) that closely approximates the translation $\tau\from N^1(\gamma) \to N^1(\gamma)$ that sends each unit normal vector to the vector obtained by parallel transport a distance one in the positive direction along $\gamma$, and then negating the normal vector. We use this permutation (and a doubling trick) to pair off the good pants, and thereby make a ``good panted surface''. Finally, we verify that any good panted surface is $\pi_1$-injective.

We can outline this approach as follows:
\begin{enumerate}
    \item 
    Define a notion of good pants (depending on parameters $\epsilon$ and $R$),
    which should be some precise version of ``the image of the $R$-perfect pants by an immersion that is $\epsilon$-close to an isometry''.
    Take one of each good pants (as in \cite{KW21}), or define some other distribution on the finite set of $(R, \epsilon)$-good pants.
    \item 
    For each good pants $\Pi$ and boundary geodesic $\gamma$,
    describe the (short) foot of $\Pi$ on $\gamma$.
    \item
    For each closed geodesic $\gamma$ in $M$ that is good
    (in the sense of being the boundary of at least one good pants),
    show that the feet of good pants are \emph{approximately invariant} by the centralizer $Z_\gamma$. We can think of this centralizer as the combination of parallel transport and the centralizer of the monodromy.
    \item 
    Conclude from this approximate invariance that there is a permutation of the feet that closely approximates the ``negate and twist by one'' translation $\tau$ on $N^1(\gamma)$.
    \item 
    Use this permutation and the ``doubling trick'' to match good pants along their common boundaries and form a good panted surface.
    \item 
    Show that every good panted surface is $\pi_1$-injective.    
\end{enumerate}

In \cite{Ham15},
the author (implicitly) follows the six-step outline above; the greatest novelty lies in a completely different approach to Step 6, using the differential geometry of rank-1 symmetric spaces (of non-compact type). 
To see why $\H^{2n}$ presents a difficulty,
we should properly understand the meaning of ``approximately invariant'' in Step 3. 
Let $\tilde{\gamma}$ be a component of the preimage of $\gamma$ in the universal cover (the symmetric space), and let $T_\gamma$ be a generator of the deck group preserving $\tilde{\gamma}$, so $\gamma = \tilde{\gamma}/T_\gamma$.
We let $Z_\gamma$ be the \emph{centralizer} of $T_\gamma$ in $\Isom(\H^{2n})$, and $Z'_\gamma$ the component of $Z_\gamma$ containing the identity. 

Here is an explicit description of $Z_\gamma'$ and $Z_\gamma$ in $\H^{2n}$.
Since $T_\gamma \in Z_\gamma' \subset Z_\gamma$, 
we can think of $Z'_\gamma$ and $Z_\gamma$ as acting on $\H^{2n}/T_\gamma$, and in fact on $N^1(\gamma)$. 
For a generic $\gamma$, the monodromy (of parallel transport) at a fiber of $N^1(\gamma)$ is conjugate to a generic element of $\SO(2n-1)$, which is orthogonally conjugate to $n-1$ blocks of 2-dimensional rotation matrices, along with a single 1 on the diagonal. The elements of $Z_\gamma'$ that fix this fiber are exactly the matrices of the same form (with arbitrary rotations). The elements of $Z_\gamma$ are the same, except that the last diagonal entry may be $\mathord-1$ rather than 1. 

This is easy to picture in the case of $n=2$: the (generic) monodromy is  a rotation around the axis corresponding to the eigenspace for the eigenvalue 1. We choose one unit vector in this axis and call it the north pole; the other is then naturally the south pole. There is then a north and south pole for each fiber, and these poles (and all latitudes) are preserved by parallel transport.
A general element of $Z_\gamma'$ shifts the fibers by some amount along $\gamma$, and then rotates each unit normal sphere by some amount around the poles, thus preserving latitude. 
For $Z_\gamma$ we add one more generator, that negates latitude and interchanges the north and south poles. 

\begin{figure}[htpb]
    \centering
    \includegraphics[width=0.4\textwidth]{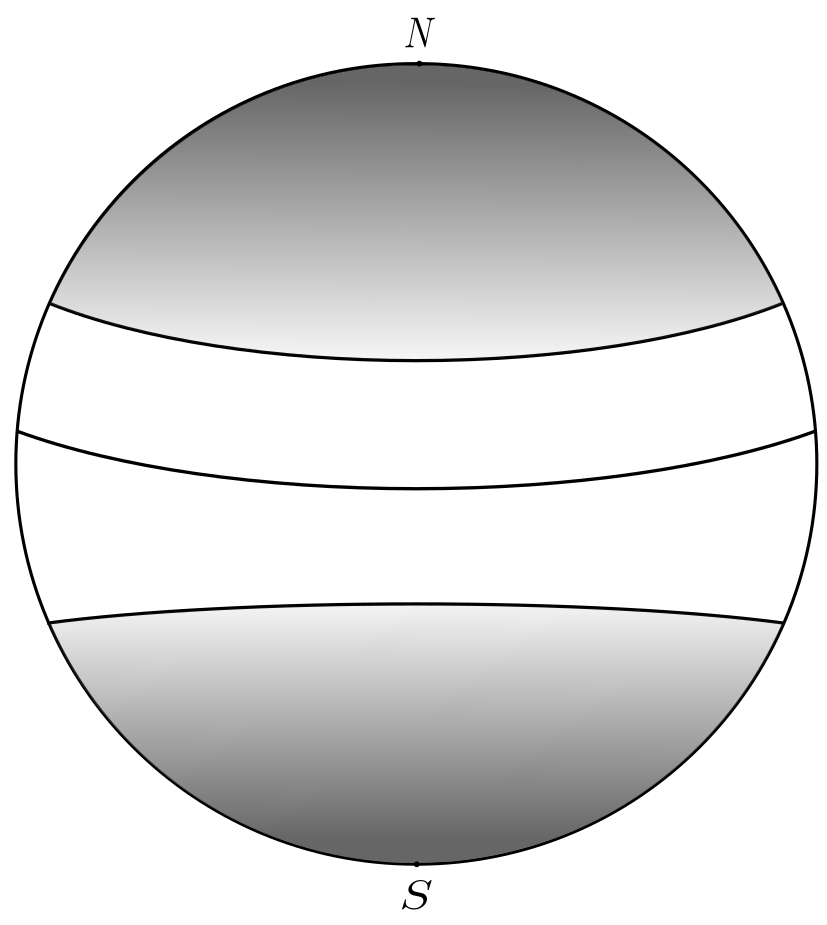}
    \caption{One possible obstruction: the ``gap" between the two components}
    \label{fig:ice_caps}
\end{figure}

In all cases, one can prove (as Step 3) that the good feet are 
\emph{approximately invariant with respect to $Z_\gamma$}. 
When a finite set of points in a compact manifold is approximately invariant with respect to a \emph{connected} group of isometries, there is a permutation of this finite set that approximates any given element of this group of isometries. When $\tau$ (defined in Step 4) lies in $Z'_\gamma$, this approximate invariance therefore implies the existence of the permutation required in Step 4. This holds in every rank one case but $\H^{2n}$.  

\begin{figure}[htpb]
    \centering
    \includegraphics[width=0.4\textwidth]{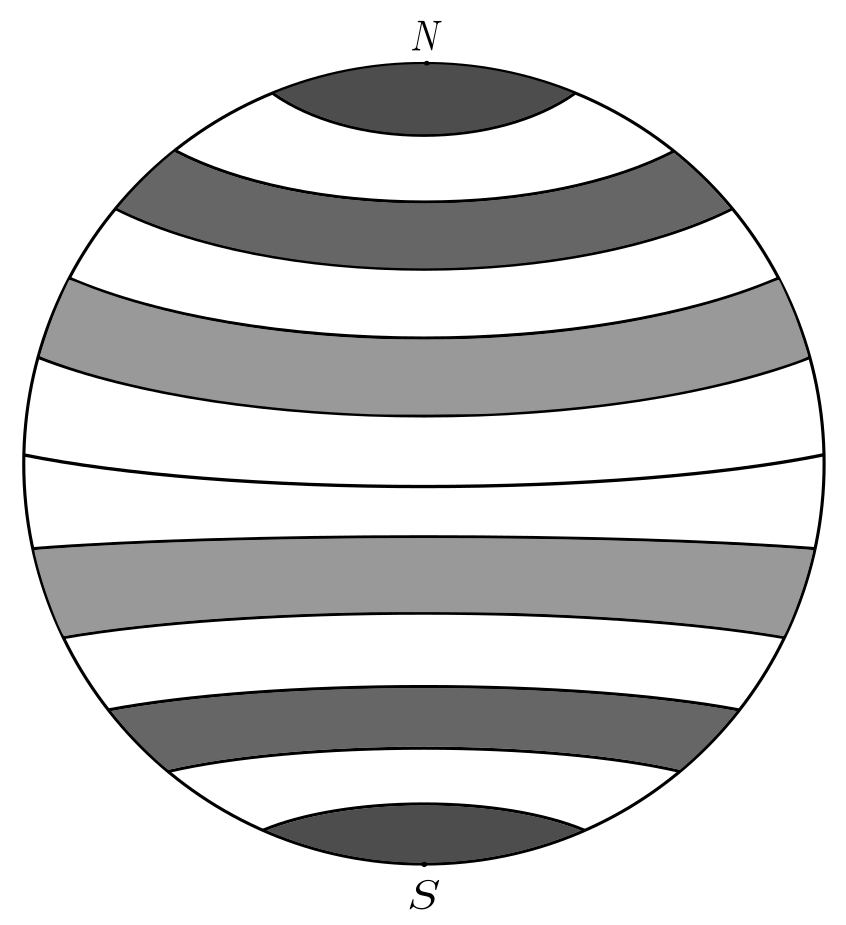}
    \caption{A worse case: the ``bands"}
    \label{fig:bands}
\end{figure}

In the case of $\H^{2n}$, this $\tau$ lies in $Z_\gamma \sm Z_\gamma'$; it is the matrix of all $\minusone$'s on the diagonal. In Figure \ref{fig:ice_caps},
in the case of $n=2$,
we illustrate a possible distribution of feet that is nearly invariant under $Z_\gamma$: in each fiber, the feet are clustered near the north and south poles, with a large gap in the lower latitudes. The problem is that there may be a few more feet in the northern ``ice cap'' than the southern one. This would prevent the existence of a bijection between the two sets of feet that would approximate the antipodal map $\tau$. (Here we present the picture in a single fiber of $N^1(\gamma)$, but the phenomenon is the same in the union of the fibers.)
Moreover, 
there could in principle be even more troublesome cases, where the feet lie in bands on the sphere, as in Figure \ref{fig:bands}.


To solve the above problem for $\H^{2n}$ with $n \ge 2$, we first observe that in Step 3, the distribution $\nu_a$ of the feet of good pants on $\gamma$ is approximated by the \emph{estimated foot measure} $\mu_a$, which can be defined analytically through third connections; the measure $\mu_a$ is \emph{exactly} invariant by $Z_\gamma$. 
(The use of third connections and this observation are inspired by \cite{KW21}.)
We then show that the Radon–Nikodym derivative of $\mu_a$ with respect to the Lebesgue measure on $N^1(\gamma)$ is bounded above and below by constants that depend only on $\epsilon$ and the manifold $M$, i.e., $\mu_a$ is \emph{quasi-uniform} with respect to the Lebesgue measure. This means that there cannot be any ``gaps'' in $\mu_a$ (or $\nu_a$) that would be an obstruction to forming the permutation that approximates $\tau$. We can then directly verify that $\nu_a$ satisfies the conditions of Hall's Marriage Theorem, and finally produce the permutation in Step 4. 


\newcommand{\anti}{\tau}
We hope,
in a sequel to this paper,
to develop a general theory of good pants in a semi-simple Lie group of non-compact type;
this would include relating two definitions of good pants, in a way similar to Theorem \ref{2.9},
and defining the feet of a good pants.
There are then, 
in general,
three possibilities of how the map $\anti$ relates to the space of feet for a closed geodesic (or flat) $\gamma$:
\begin{enumerate}
    \item 
    $\anti$ lies in the identity component $Z_\gamma'$ of the centralizer for $T_\gamma$;
    \item 
    $\anti$ lies in the identity component of the stabilizer of $\hat\gamma$ (the cover of $\gamma$), but not in $Z_\gamma'$;
    \item 
    $\anti$ does not lie in the identity component of the stabilizer of $\gamma$.
\end{enumerate}
The first case is the only one treated so far, in \cite{Ham15} and \cite{KLM24} (and \cite{KM12b} and \cite{KW21}).
The third case would definitely require the correction theory of \cite{KM15}.
This paper is the first paper treating the second case;
the method developed here, 
using quasi-uniformity of the foot measure,
should be applicable to the entirety of the second case,
and may also be a necessary preliminary to the third.

While the \emph{problem} arises only for even-dimensional cases,
the \emph{solution} works perfectly well for any dimension at least 4;
we therefore denote the dimension by $n$ from now on and treat all $n$ on an equal basis, 
no longer concerning ourselves with whether $n$ is odd or even.

Section \ref{sec:geometry_of_pants} introduces preliminaries on the geometry of pants, especially in $n$-dimensional hyperbolic space and manifolds. 
After recalling properties of $\SO(n,1)$, and basic definitions for pants, 
we define the notion of \emph{cuff-good} pants, 
which are pants whose cuffs are good curves (as in \cite{KW21}). 
Using the \emph{Steiner graph} (which we define) for pants, 
we then show that cuff-good pants are either good (in the sense of  \cite{Ham15} and \cite{KM12b}, where we generate pants from connections between tripods), or \emph{bad}, in that they are close to a geodesic immersion of a one-holed torus. 
This then allows us to determine when a pair of pants is good from one cuff and a \emph{third connection} that generates the other two cuffs.

Section \ref{sec:counting} contains the core estimate of this paper, where we approximate the distribution of cuff-good pants through counting the third connections. Given a closed geodesic $\gamma$, for any framed third connection, we compute the monodromy of the other two cuffs of the pants using matrix multiplication, and then derive an (almost) necessary and sufficient condition for generating cuff-good pants. Thus we can determine the \emph{good region} in the space of third connections, which is a fiber bundle over the unit normal bundle $N^1(\gamma)$ of $\gamma$. The projection of the good region yields a measure on $N^1(\gamma)$, which is called the \emph{estimated foot measure}. By evaluating the volume of the good region in each fiber, we show that the estimated foot measure is bounded above and below by a multiple of the Lebesgue measure; this rules out the ``ice cap'' or ``band'' obstruction described above.

In Section \ref{sec:final-proof}, we first show that our matching pattern satisfies the incompressibility criterion established by Hamenst{\"a}dt\cite{Ham15}. We then apply Hall's Marriage Theorem (using our measure bounds and estimates of the Cheeger constant) to produce the permutation on the set of all good pants, which is used to properly match the good pants and prove the main result of this paper, the existence of $\pi_1$-injective immersed surfaces in closed hyperbolic $n$-manifolds.

\textbf{Acknowledgments.} This work arose out of a conversation in 2014, in which Ursula Hamenst{\"a}dt explained to the first author the problem that can arise when trying to build a surface out of pants in a quotient of $\H^4$. 

The second author would like to thank Fernando Al Assal for inviting him to visit MPI MIS in 2023. In a conversation during this visit, Al Assal believed that pants with good cuffs should be good. Although the claim is now known to be not entirely true, it inspired the second author to think about the monodromy of cuffs. Finally, the first and second authors managed to show that a pants in $\H^{2n}$ with good cuffs is either good or bad. The second author would also like to thank Ursula Hamenst{\"a}dt and Yi Liu for helpful conversations.

The authors owe thanks to Hongbin Sun for detailed corrections and suggestions for this paper.

\section{Some geometry of good pants}\label{sec:geometry_of_pants}

\subsection{The structure of \texorpdfstring{$\SO^{+}(n,1)$ and $\mathfrak{so}(n,1)$}{}}

Let $\mathcal{G}=\SO^{+}(n,1)$, which is the identity component of $\SO(n,1)$ and acts on $\H^{n}$ by orientation preserving isometries. We choose $p_0=(1, 0,\cdots,0)\in\H^{n}$ (in the upper half space model) to be the base point, and let $\mathbf{u}_i=(0,\cdots,1, \cdots, 0)$ be the unit tangent vector at $p_0$ with a 1 at the $i^{\rm{th}}$ place. We let $F_0 = (u_1, \ldots, u_{n})$ be the base frame at $p_0$. The left action of $\mathcal{G}$ on $\H^{n}$ as isometries induces a left action of $\mathcal{G}$ on the frame bundle $\mathcal{F}(\H^{n})$. 
We will think of $\SO(n)$ as the subgroup of $\SO(n, 1)$ that fixes $p_0$, and acts on $T_{p_0}(\H^{n}) = \R^{n}$ in the usual way.

Now we define the right action of $\mathcal{G}$ on $\mathcal{F}(\H^{n})$ based on the left action, generalizing the definition in Section 2 of \cite{Rao25} for dimension 2 and 3. Let $\Psi_0=\langle p_0,F_0\rangle\in\mathcal{F}(\H^{n})$ be the base frame. Then the right action of $\mathcal{G}$ on $\mathcal{F}(\H^{n})$ is defined as follows:
\begin{enumerate}[(1)]
    \item For any $g\in\mathcal{G}$, $\Psi_0\cdot g=g\cdot\Psi_0$.

    \item For any $\Psi\in\mathcal{F}(\H^{n})$ and $g,h\in\mathcal{G}$, $g\cdot(\Psi\cdot h)=(g\cdot\Psi)\cdot h$.
\end{enumerate}
Since the left action of $\mathcal{G}$ on $\mathcal{F}(\H^{n})$ is transitive and free, then for any $\Psi$, there exists a unique $g\in\mathcal{G}$ such that $\Psi=g\cdot\Psi_0$. Thus for any $h\in\mathcal{G}$, 
\begin{equation*}
    \Psi\cdot h=(g\cdot\Psi_0)\cdot h=g\cdot(\Psi_0\cdot h),
\end{equation*}
which means the right action on $\Psi_0$ will tell us the right action on any other frames.
Thus $\mathcal{G}$ acts on the right by instructions, such as frame flows and rotations.
In particular,
every matrix in $\SO(n) < \mathcal{G}$ gives an instruction that creates a new frame in terms of the old frame without changing its base point. For example, $\mathbf{diag}(-1, -1, 1, \ldots, 1)$ negate the first two vectors in the frame, and the matrix  
$\mathbf{diag}\left(\begin{pmatrix}
0 & -1 \\
1 & 0
\end{pmatrix},1,\ldots,1\right)$
rotates the first two vectors of the frame counter-clockwise by $\pi/2$, \emph{in the plane of those first two vectors}. 
We say that the matrix acts on the frame by ``rewriting'', 
in the sense that columns of the matrix determine how the vectors of the new frame are written as linear combinations of vectors of the old frame.  

Let $\mathcal{K}$ be the stabilizer of $p_0$ (for $\mathcal{G}$ acting on the left); then $\mathcal{K}$, which we have identified with $\SO(n)$, is a maximal compact subgroup of $\mathcal{G}$. 
The stabilizer $\mathcal{M}$ of $\mathbf{u}_0$ is a subgroup of $\mathcal{K}$ which is isomorphic to $\SO(n-1)$. 
We've specified how $\mathcal{K}$ and $\mathcal{M}$ act on the \emph{left}; 
we have also observed in the previous paragraph that $\mathcal{K}$ acts on the right by rewriting the frame, and $\mathcal{M}$ acts by rewriting but preserving the first vector in the frame. 
Hence $\mathcal{G}/\mathcal{K}$ (where we take the quotient by the right action) can be identified with $\H^{n}$, $\mathcal{G}/\mathcal{M}$ can be identified with the unit tangent bundle $T^1(\H^{n})$ of $\H^{n}$, and $\mathcal{G}$ itself can be identified with the oriented frame bundle $\mathcal{F}(\H^{n})$ of $\H^{n}$. We always denote by $\mathbf{e}$ the identity element of a group, and $[\cdot]$ the conjugacy class of a group element.

There is a one parameter group $\mathcal{A}=\{a_t:t\in\R\}$ such that the right action of $a_t$ on $\mathcal{G}$ is the unit-speed frame flow, and such that $\mathcal{A}$ commutes with $\mathcal{M}$. Let $\mathcal{B}=\mathcal{A}\mathcal{M}=\mathcal{M}\mathcal{A}<\mathcal{G}$. We can define the unstable and stable horospherical subgroups $\mathcal{N}^+$ and $\mathcal{N}^-$ of $\mathcal{G}$, respectively, as follows:
\begin{equation*}
    \mathcal{N}^{\pm}=\{g\in\mathcal{G}: a_t g a_{-t}\to \mathbf{e}, \mathrm{as}\ t\to\pm\infty\}
\end{equation*}
Let $\mathfrak{a},\mathfrak{m},\mathfrak{n}^{\pm}$ be the corresponding Lie subalgebra in $\mathfrak{so}(n,1)$ of $\mathcal{A},\mathcal{M},\mathcal{N}^{\pm}$, respectively. Thus
\begin{equation*}
    \mathfrak{so}(n,1)=\mathfrak{a}\oplus\mathfrak{m}\oplus\mathfrak{n}^+\oplus
    \mathfrak{n}^-,
\end{equation*}
as real vector spaces.

By hyperbolic geometry, we have the following facts:
\begin{enumerate}
    
    \item For any $X^{+}\in\mathfrak{n}^+$ and $t>0$,
    \begin{equation*} 
        a_t\exp(X^+)=\exp(e^{-t}X^+)a_t.
    \end{equation*}

    \item For any $X^{-}\in\mathfrak{n}^-$ and $t>0$,
    \begin{equation*}
        \exp(X^-)a_t=a_t\exp(e^{-t}X^-).
    \end{equation*}
\end{enumerate}

\begin{lem}\label{2.1}
    There exist $C,\epsilon_0>0$ such that for any $\epsilon \in (0, \epsilon_0)$: Suppose that $u\in\mathcal{G}$ is $\epsilon$-close to $\mathbf{e}$, then $u$ can be uniquely written as 
    \begin{equation*}
        u=n^+ b n^-,
    \end{equation*}
    for some $n^\pm\in\mathcal{N}^\pm$ and $b\in\mathcal{B}$. Moreover, $n^\pm$ and $b$ are $C\epsilon$-close to $\mathbf{e}$. 
\end{lem}

For $X, Y \in \mathfrak{g}$, we let $B(X, Y)$ be the Killing form. As usual we can write $\mathfrak{g} = \mathfrak{k} \oplus \mathfrak{p}$, where $\mathfrak{p}$ is the subspace of infinitesimal transvections at $p_0$. For any $X = X_{\mathfrak{k}} + X_{\mathfrak{p}} \in \mathfrak{k} + \mathfrak{p}$, we let $\theta(X) = X_{\mathfrak{k}} - X_{\mathfrak{p}}$ be the Cartan involution; we might then define the inner product $\pair{X}{Y} := -B(X, \theta(Y))$. 
The infinitesimal rotation $\Theta := \mathbf{diag}\left(\begin{pmatrix}
0 & -1 \\
1 & 0
\end{pmatrix},1,\ldots,1\right) \in \mathfrak{k}$ has $B(\Theta, \Theta) = -2(n-2)$, but we want it to have norm 1, so we actually let $\pair{X}{Y} := -B(X, \theta(Y))/(2n-4)$. 
This then induces a left-invariant Riemannian metric on $\mathcal{G}$;
it is also right-invariant by $\mathcal{K}$,
and therefore restricts to a bi-invariant metric on $\mathcal{K}$.
For $v, w \in S^{n-1}$, 
let $d(v, w) \in [0, \pi]$ be the spherical distance. 
We have the following estimates, for any $X, Y \in \SO(n)$ and $v \in S^{n-1}$, which are left to the reader:
\begin{enumerate}
    \item $d(Xv, Yv) \le d(X, Y)$, 
    \item $d(X, \mathbf{e}) = d(X v, v)$ when $X$ is a rotation of a 2-plane $V$ (and the identity on $V^\perp$) and $v \in V$. 

\end{enumerate}
We let $\|\cdot \|$ be the norm corresponding to the inner product.

\begin{lem}\label{2.2}
    For any $X^+\in\mathfrak{n}^+$, $t>0$ and $m\in\mathcal{M}$, there exists $X^+_1\in\mathfrak{n}^+$ such that
    \begin{equation*}
        \|X^+_1\|=e^{-t}\|X^+\|,
    \end{equation*}
    and
    \begin{equation*}
        a_t m\exp(X^+)=\exp(X^+_1) a_t m.
    \end{equation*}
\end{lem}

\begin{proof}
    Let $X^+_1=\Ad_m (e^{-t}X^+)$, then $\|X^+_1\|=e^{-t}\|X^+\|$, and we have
    \begin{equation*}
        \begin{split}
            a_t m \exp(X^+) (a_tm)^{-1}&=m a_t \exp(X^+) a_t^{-1} m^{-1}\\
            &=m \exp(e^{-t}X^+) m^{-1}\\
            &=\exp(\Ad_m(e^{-t}X^+))=\exp(X^+_1).\qedhere
        \end{split}
    \end{equation*}
\end{proof}

Similarly, we have another lemma.

\begin{lem}\label{2.3}
    For any $X^-\in\mathfrak{n}^-$, $t>0$ and $m\in\mathcal{M}$, there exists $X^-_1$ such that
    \begin{equation*}
        \|X^-_1\|=e^{-t}\|X^-\|,
    \end{equation*}
    and
    \begin{equation*}
        \exp(X^-)a_t m=a_t m\exp(X^-_1).
    \end{equation*}
\end{lem}

Suppose that $g \in \mathcal{G}$ preserves a geodesic $\tilde{\gamma}$ through $(p_0, u_1)$; then $g \in \mathcal{B}$ (as a left action). Of course we then have $g \cdot \left< p_0, F_0 \right> = \left< p_0, F_0 \right> \cdot g$ by the definition of the right action. Writing $g = a_t m$, we see that the monodromy of $\gamma:= \tilde{\gamma}/g$, as a left action on $T_{p_0} H^{n}$, is given by $m$. Moreover the image of $F_0$ by this monodromy is $m \cdot F_0 = F_0 \cdot m$.

More generally, suppose that $g \in \mathcal{G}$ (as an isometry) translates along (and possibly rotates around) a geodesic $\tilde{\gamma}$ with quotient $\gamma$. Suppose that $\tilde{p} \in \tilde{\gamma}$, and $\tilde{F}$ is a frame at $\tilde{p}$, with first vector tangent to $\tilde{\gamma}$ in the direction of translation by $g$. 
Let $p$ and $F$ be the image of $\tilde{p}$ and $\tilde{F}$ on $\gamma$.
Then the monodromy of $\gamma$ maps $F$ to another frame, which must have the form $F \cdot m$, for $m \in \mathcal{M}$. Letting $t$ be the length of $\gamma$, we claim that $g$ is conjugate to $a_t m$.
To see this, we take $h\in \mathcal{G}$ such that $h \cdot (p_0, F_0) = (\tilde{p}, \tilde{F})$; then $$h^{-1}gh \cdot (p_0, F_0) = 
h^{-1}g \cdot(\tilde{p}, \tilde{F}) =h^{-1} \cdot (\tilde{p}, \tilde{F}) \cdot a_t m = 
(p_0, F_0) \cdot a_t m = a_t m \cdot (p_0, F_0),$$
so $h^{-1} g h = a_t m$.

In the following two lemmas,
we consider certain classes of words in $\mathcal{G}$ that must preserve a geodesic,
and we provide an estimate of this length and monodromy for the quotient.
The first one tells us that a small perturbation can be absorbed into long frame flows and rotations after conjugacy, and the second generalizes the first to multiple elements and perturbations.

\begin{lem}\label{2.4}
    There exist $\epsilon_0,R,K>0$ such that when $\epsilon_0>\epsilon>0$, the following holds: Suppose that $u\in\mathcal{G}$ is $\epsilon$-close to $\mathbf{e}$, $t>R$ and $m\in\mathcal{M}$. Then there exists $t'>0$ with $|t-t'|<K\epsilon$ and $m'\in\mathcal{M}$ that is $K\epsilon$-close to $m$, such that
    $a_tmu$ is conjugate to $a_{t'}m'$.
\end{lem}

\begin{proof}
    Suppose $\epsilon_0$ and $C$ are the constants in Lemma \ref{2.1}, and $R>10000$ large enough, which will be determined later. Next we will find $K$ satisfying the requirements.

    By Lemma \ref{2.1}, there exist $n^\pm_1\in\mathcal{N}^\pm$, $a_{t_1}\in\mathcal{A}$ and $m_1\in\mathcal{M}$ such that 
    \begin{equation*}
        u=n^+_1 a_{\hat{t}_1} \hat{m}_1 n^-_1
    \end{equation*}
    and $n^\pm_1,\hat{m}_1,a_{\hat{t}_1}$ are all $C\epsilon$-close to $\mathbf{e}$. Then by Lemma \ref{2.2} and Lemma \ref{2.3}, we have
    \begin{equation*}
        \begin{split}
            [a_tmu]&=[a_tmn^+_1 a_{\hat{t}_1} \hat{m}_1 n^-_1]=[a_{t/2}mn^+_1 a_{\hat{t}_1} \hat{m}_1 n^-_1a_{t/2}]\\
            &=[\hat{n}_1^+a_{t/2}ma_{\hat{t}_1} \hat{m}_1 a_{t/2}\hat{n
            }_1^-]\\
            &=[a_{t+\hat{t}_1}m\hat{m}_1\hat{n}_1^-\hat{n}_1^+],
        \end{split}        
    \end{equation*}
    where $\hat{n}_1^-$ and $\hat{n}_1^+$ are $e^{-t}C\epsilon$-close to $\mathbf{e}$.
    Now we let $t_1=t+\hat{t}_1$, $m_1=m\hat{m}_1$ and $u_1=\hat{n}_1^-\hat{n}_1^+$. Thus $|t-t_1|<C\epsilon$, $m_1$ is $C\epsilon$-close to $m$, $u_1$ is $2e^{-t}C\epsilon$-close to $\mathbf{e}$, and
    \begin{equation*}
        [a_tmu]=[a_{t_1}m_1u_1].
    \end{equation*}
    Let $K=2C$, and $R$ large enough such that $2e^{-R/2}C<1/2$ and $t_i>R/2$ for all $i\in\Z_+$. Thus we can have sequences $\{t_i\}_{i=1}^\infty\subset \R_+$, $\{m_i\}_{i=1}^\infty\subset \mathcal{M}$ and $\{u_i\}_{i=1}^\infty\subset\mathcal{G}$ such that for any $i\in\Z_+$, $|t_i-t_{i+1}|<2^{1-i}C\epsilon$, $m_{i+1}$ is $2^{1-i}C\epsilon$-close to $m_i$, $u_{i}$ is $2^{-i}\epsilon$-close to $\mathbf{e}$, and
    \begin{equation*}
        [a_{t_i}m_iu_i]=[a_{t_{i+1}}m_{i+1}u_{i+1}].
    \end{equation*}
    Therefore the limits of these three sequences exist, and moreover, the limit of $\{u_i\}$ is $\mathbf{e}$. Let $t_i\to t'$ and $m_i\to m'$, as $i\to\infty$. Hence $[a_tmu]=[a_{t'}m']$, $|t-t'|<2C\epsilon$, and $m'$ is $2C\epsilon$-close to $m$, which completes the proof.
\end{proof}
We now estimate the length and monodromy in a more complicated case. 
\begin{lem}\label{8-word}
    There exist $R,\epsilon_0>0$ and $K>1$ such that for any $\epsilon_0>\epsilon>0$, the following holds: Suppose $u_1,u_2,v_1,v_2\in\mathcal{G}$ are $\epsilon$-close to $\mathbf{e}$, $t_1,t_2>R$, and $m_1,m_2\in\mathcal{M}$. Then there exists $t\in\R$ and $m\in\mathcal{M}$ such that:
    \begin{enumerate}
        \item $|t-(t_1+t_2)|<K\epsilon$;
        \item $m$ is $K\epsilon$-close to $m_1m_2$;
        \item $[a_{t_1} u_1 m_1 v_1 a_{t_2} u_2 m_2 v_2]=[a_t m]$.
    \end{enumerate}

\end{lem}

\begin{proof}
    Let $R$ be the same constant $R$ in the previous lemma, $\epsilon_1$ and $C$ be the constants in Lemma \ref{2.1} (which is the same as in Lemma \ref{2.4}), and $\epsilon_0<\epsilon_1$ which will be determined later. Then as in the proof of Lemma \ref{2.4}, we have the decomposition $u_i=n^+_i a_{\hat{t}_i} \hat{m}_i n^-_i$ for $i=1,2$. Then by Lemma \ref{2.2} and Lemma \ref{2.3}, we have
    \begin{equation*}
        \begin{split}
            [a_{t_1} u_1 m_1 v_1 a_{t_2} u_2 m_2 v_2]&=[a_{t_1}n^+_1 a_{\hat{t}_1} \hat{m}_1 n^-_1 m_1 v_1 a_{t_2} n^+_2 a_{\hat{t}_2} \hat{m}_2 n^-_2 m_2 v_2]\\
            &=[\hat{n}^+_1 a_{t_1} a_{\hat{t}_1} \hat{m}_1 m_1 \hat{n}^-_1 v_1 \hat{n}^+_2 a_{t_2} a_{\hat{t}_2} \hat{m}_2 m_2 \hat{n}^-_2 v_2]\\
            &=[a_{t_1+\hat{t}_1}\hat{m}_1 m_1 u_3 a_{t_2+\hat{t}_2}\hat{m}_2 m_2 u_4],
        \end{split}
    \end{equation*}
    where $\hat{n}_i^+$'s are $Ce^{-R}\epsilon$-close to $\mathbf{e}$, $\hat{n}_i^-$'s are $C\epsilon$-close to $\mathbf{e}$, and
    \begin{equation*}
        u_3=\hat{n}^-_1 v_1 \hat{n}^+_2, u_4=\hat{n}^-_2 v_2 \hat{n}^+_1
    \end{equation*}
    are $(1+C+Ce^{-R})\epsilon$-close to $\mathbf{e}$. Then we can repeat the above process to decompose $u_3$ as $u_3=n^+_3 a_{\hat{t}_3} \hat{m}_3 n^-_3$, and we will have
    \begin{equation*}
        [a_{t_1} u_1 m_1 v_1 a_{t_2} u_2 m_2 v_2]=[a_{t_1+\hat{t}_1+\hat{t}_3+t_2+\hat{t}_2} \hat{m}_1 m_1 \hat{m}_3 \hat{m}_2 m_2 \hat{n}_3^- u_4 \hat{n}^+_3],
    \end{equation*}
    where $\hat{n}_3^\pm$ are $Ce^{-R}(1+C+Ce^{-R})\epsilon$-close to $\mathbf{e}$. Now let $t_5=t_1+\hat{t}_1+\hat{t}_3+t_2+\hat{t}_2$, $m_5=\hat{m}_1 m_1 \hat{m}_3 \hat{m}_2 m_2$ and $u_5=\hat{n}_3^- u_4 \hat{n}^+_3$, and we know the distance between $u_5$ and $\mathbf{e}$ is at most
    \begin{equation*}
        (1+C+Ce^{-R})\epsilon+2Ce^{-R}(1+C+Ce^{-R})\epsilon.
    \end{equation*}
    Thus we choose $\epsilon_0<\epsilon_1$ such that 
    \begin{equation*}
        (1+C+Ce^{-R})\epsilon_0+2Ce^{-R}(1+C+Ce^{-R})\epsilon_0=\epsilon_1,
    \end{equation*}
    and then we have $u_5$ is $\epsilon_1$-close to $\mathbf{e}$. By Lemma \ref{2.4}, we can find $t\in\R_+,m\in\mathcal{M}$ that are $2\epsilon_1$-close to $t_5,m_5$, respectively, such that $[a_{t_5}m_5u_5]=[a_tm]$.
\end{proof}

\subsection{Preliminary definitions}
We want to introduce some important geometric definitions related with pants and surfaces in hyperbolic n-manifolds. Some are recalled from \cite{KM12b}, \cite{Ham15} and \cite{Rao25}.

\subsubsection{Monodromy, good curves and normal bundles}
Suppose $M$ is a hyperbolic n-manifold and $\gamma\subset M$ is a closed geodesic. Then $\gamma$ is the quotient of a geodesic in $\H^{n}$ under some isometry of $\H^{n}$. Such an isometry is uniquely determined up to conjugacy by its translation length, which is the length of $\gamma$, and its rotation part, which is an element in $\SO(n-1)$. The length of $\gamma$, denoted by $l(\gamma)$, is called $\epsilon$-\textit{close} to $R$, if
\begin{equation*}
    R-\epsilon<l(\gamma)<R+\epsilon,
\end{equation*}
for some $R,\epsilon>0$.

The rotation part is called the \emph{monodromy} of $\gamma$. The monodromy of $\gamma$ is called $\epsilon$-\textit{close} to $\mathbf{e}$, if distance between the monodromy and $\mathbf{e}$ is bounded by $\epsilon$ in $\SO(n-1)$. We then say that a closed geodesic $\gamma$ is $(R, \epsilon)$-\emph{good} if it has real length $2\epsilon$-close to $2R$ and monodromy $\epsilon$-close to $\mathbf{e}$.

\begin{rem}
    For $g\in \SO(m)$, we can regard that $g$ is $\epsilon$-close to $\mathbf{e}$ for some $\epsilon>0$, if the angle between $v$ and $gv$ is no greater than $\epsilon$ for any unit vector $v$. For $g,g'\in\SO(m)$, we say $g$ is $\epsilon$-close to $g'$, if $g^{-1}g'$ is $\epsilon$-close to $\mathbf{e}$.
\end{rem}

Suppose $\gamma\colon U\to M$ is a unit speed geodesic in $M$, where $U$ is a Riemannian 1-manifold. The unit normal bundle of $\gamma$ is defined as
\begin{equation*}
    N^1(\gamma)=\{(u,v):u\in U,v\in T_{\gamma(u)}M,\langle v,\gamma'(u)\rangle=0,||v||=1\}.
\end{equation*}
If $\gamma$ is not closed, $N^1(\gamma)$ is isometric to $[0,l(\gamma)]\times S^{n-2}$. Otherwise $N^1(\gamma)$ is the quotient of $\R\times S^{n-2}$ by the translation along $\gamma$. In both cases, $N^1(\gamma)$ can be regarded as a metric space where the metric comes from the $L^2$-metric on $\R\times S^{n-2}$. Moreover, we define the normal frame bundle of $\gamma$ as follows:
\begin{equation}
    \mathcal{F}^\perp(\gamma)=\{(u,E):u\in U,(\gamma'(u),E)\ \text{is\ a\ positive-oriented\ frame\ of\ } T_uM\}.
\end{equation}

\subsubsection{Orthogeodesics, connections and feet}\label{sub4.3}
Suppose there are two oriented closed geodesics $\alpha_i$, $i=1,2$, in $M$. Suppose that $\eta\colon[0,l]\to M$ is a unit speed geodesic segment such that
\begin{itemize}
\item $\eta(0)\in\alpha_1$, $\eta(l)\in\alpha_2$, and
\item $\eta(0,l)$ meets with $\alpha_i$'s orthogonally, $i=1,2$.
\end{itemize}
We then call $\eta$ an \emph{orthogeodesic} or a \emph{connection} between $\alpha_1$ and $\alpha_2$. In the case where $\alpha_1$ and $\alpha_2$ coincide, we will call $\eta$ a \emph{third connection}. 

If $l > 0$, we let $\mathrm{foot}_{\alpha_1}(\eta)={\eta}'(0)\in N^1(\alpha_1)$ and $\mathrm{foot}_{\alpha_2}(\eta)=-{\eta}'(l)\in N^1(\alpha_2)$ be the \emph{feet} of $\eta$ on $\alpha_0$ and $\alpha_1$, respectively. 
The base point of a foot is called a \emph{footpoint}.
(If $l = 0$,
the footpoints are defined, 
but not the feet.
This case will not occur in the good pants that we are studying.)

Similarly, we can define the orthogeodesic and feet for two  geodesics in $\H^{n}$ at positive distance from each other;  the orthogeodesic between them is always unique.

\subsubsection{Pants and feet}

Suppose $\Pi_0$ is a topological pair of pants considered as a manifold with boundary. We say a \emph{pants} in $M$ is an injective homomorphism $\rho\colon\pi_1(\Pi_0)\to \pi_1(M)$, up to conjugacy, which determines and is determined by a continuous map $f\colon \Pi\to M$, up to homotopy. So we can also call $f$ or $f(\Pi_0)$ a pants. Let $C_i$'s be the boundary components of $\Pi_0$, and we orient $C_i$ such that $\Pi_0$ is on the left of $C_i$, for $i\in\Z/3\Z$. Then there is a unique oriented closed geodesic $\gamma_i$ in $M$ that is freely homotopic to $f(C_i)$. We can homotope $f$ such that $f$ maps $C_i$ to $\gamma_i$, and we call such $f$ a \emph{nice} pants. From now on, when talking about pants in $M$, we will always assume that they are nice.

For any given simple non-separating arc $\alpha_i$ in $\Pi_0$ connecting $C_{i-1}$ and $C_{i+1}$, we then can homotope $f$ such that $f$ maps $\alpha_i$ to an orthogeodesic $\eta_i$ between $\gamma_{i-1}$ and $\gamma_{i+1}$. Here $\eta_i$ does not depend on the choice of $\alpha_i$ and $f$ within its homotopy class. We call $\gamma_i$'s and $\eta_i$'s the \emph{cuffs} and the \emph{short\ orthogeodesic} of the pants $f(\Pi_0)$, respectively. The feet of $\eta_{i-1}$ and $\eta_{i+1}$ on $\gamma_i$ are called the \emph{short feet} of $f$ on $\gamma_i$. When we talk about a cuff $\gamma$, the information of the associated pants is also carried by $\gamma$. 

Moreover, for any simple arc $\beta_i$ in $\Pi_0$ connecting $C_i$ and itself such that $C_{i-1}$ and $C_{i+1}$ are in different component of $\Pi_0-\beta_i$ (which always has two connected components), we can homotope $f$ such that $f$ sends $\beta_i$ to an orthogeodesic $\zeta_i$ from $\gamma_i$ to itself. $\zeta_i$ is independent of the choice of $\beta_i$ and $f$, and we also call $\zeta_i$ is the third connection for $\gamma_i$ in $\Pi_0$. The two feet of $\zeta_i$ on $\gamma_i$ are called the \emph{long feet} of $f$ on $\gamma_i$.

\subsubsection{Tripods and framed tripods}
We define a \emph{tripod} in $T\H^{n}$ or $TM$ to be an ordered triple $(v_0,v_1,v_2)$ of unit tangent vectors at the same base point $x$ such that the three vectors mutually enclose an angle of $2\pi/3$. The vectors of a tripod are contained in a 2-plane $V$ in $T_x\H^{n}$ or $T_x M$. The orientation of $V$ is determined by the cyclic order of the triple. Given a tripod, we then let $v_i^\perp$ be the unit vector in $V$ that is orthogonal to $v_i$, and such that $(v_i, v_i^\perp)$ is a positive frame for $V$. This implies that $v_i^\perp$ is at an angle of $\pi/6$ to $v_{i+1}$. 

A \emph{framed tripod} in $T\H^{n}$ or $TM$ is a pair $((v_0,v_1,v_2),E)$ where $(v_0,v_1,v_2)$ is a tripod contained in a plane $V$ and $E$ is a positive orthonormal frame in the orthogonal complement of $V$. We then can define how two framed tripods are related to each other.

The following definition is combined from Section 4 of \cite{KM12b} and Section 4 of \cite{Ham15}, which plays an essential role in the definition of good pants.

\begin{defn} \label{good-conn}
    For $R>0,\epsilon>0$, we say two framed tripods $((v_0,v_1,v_2),E)$ at $p$ and $((u_0,u_1,u_2),F)$ at $q$ in $TM$ are $(R,\epsilon)$-\emph{well}-\emph{connected}, if there are geodesic segments $\gamma_i$, $i\in\Z/3\Z$, connecting $p$ and $q$, such that:
    \begin{enumerate}[(1)]
        \item The real length $l(\gamma_i)$ of $\gamma_i$ satisfies that
        \begin{equation}
            R+\log\frac43-\epsilon < l(\gamma_i) < R+\log\frac43+\epsilon;
        \end{equation}

        \item The angle between $v_i$ and $\gamma'_i$ at $p$ is bounded by $e^{-R/4}$;

        \item The angle between $u_{-i}$ and $-\gamma'_i$ at $q$ is bounded by $e^{-R/4}$;

        \item Let $E_i=(v_i,v_i^{\perp},E)$ (so $E_i$ is a positive frame at $p$), and the same for $F_i$ at $q$. Let $\hat{E_i}$ be the frame from the parallel transport of $E_i$ along $\gamma_i$ to $q$ by negating the first two vectors. Then the element in $\SO(n)$ which transforms $\hat{E_i}$ to $F_i$ is $\epsilon$-close to $\mathbf{e}$.
    \end{enumerate}
    The geodesic segments $\gamma_i$ are called \emph{good connections} of the tripods, and the transformation defined in (4) is called the monodromy of $\gamma_i$.
\end{defn}

We also want to single out another family of pairs of tripods, where the monodromy of each connection behaves more like an involution, rather than being close to $\mathbf{e}$.

\begin{defn} \label{bad-conn}
    For $R>0,\epsilon>0>0$, we say two framed tripods $((v_0,v_1,v_2),E)$ at $p$ and $((u_0,u_1,u_2),F)$ at $q$ in $TM$ are $(R,\epsilon)$-\emph{badly}-\emph{connected}, if there are geodesic segments $\gamma_i$, $i\in\Z/3\Z$, connecting $p$ and $q$, such that items (1)--(3) in Definition \ref{good-conn} hold, and:
    \begin{enumerate}[(1)]



        \item[(4')] 
        With $\hat E_i$ and $F_i$ defined as in (4) of Definition \ref{good-conn},
        the element in $\SO(n)$ which transforms $\hat{E_i}$ to $F_i$ is $\epsilon$-close to $\mathbf{diag}(1,-1,-1,1,\ldots,1)$.
    \end{enumerate}
    We then say that the $\gamma_i$ are \emph{bad connections}.
\end{defn}


\begin{figure}[htpb]
    \centering
    \includegraphics[width=0.6\textwidth]{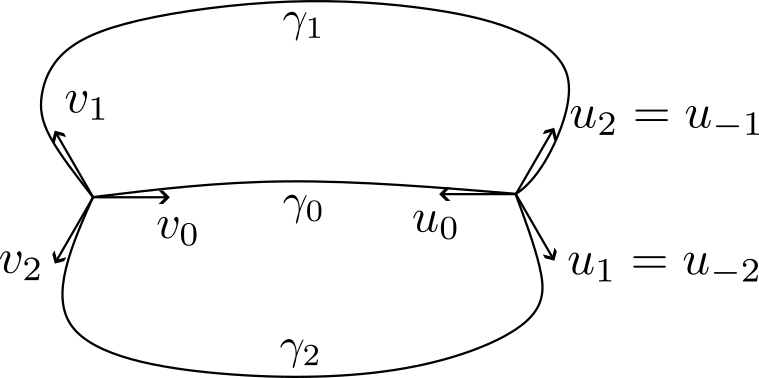}
    \caption{A pair of well-connected tripods}
    \label{fig:well-connected}
\end{figure}

\begin{figure}[htpb]
    \centering
    \includegraphics[width=0.8\textwidth]{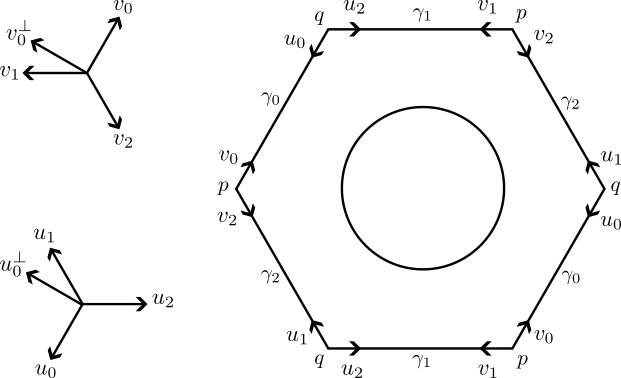}
    \caption{A pair of badly-connected tripods}
    \label{fig:badly-connected}
\end{figure}

When the tripods are well-connected, the connections generate a pants group that is close to a Fuchsian group for a pants with three cuffs of length $2R$, as showed in Figure \ref{fig:well-connected}. When the tripods are badly-connected, we obtain a group that is close to a Fuschsian group for a one-holed torus with $\Z/6\Z$ rotational symmetry and a single cuff of length $6R$. In this latter case, the orientation of $(u_0, u_1, u_2)$ is negative \emph{with respect to the plane of the Fuchsian group}, so each monodromy is orientation-reversing from the plane of the $v_i$'s to that of the $u_i$'s. Since the $\hat E_i$ and the $F_i$ are positive frames, we must also reverse a vector in the parallel transport of $E$, and our convention is to reverse the first vector (which is then the third in $\hat E_i$). This case is illustrated in Figure \ref{fig:badly-connected}, where the one-holed torus is obtained by identifying the opposite sides of the hexagon.

\subsubsection{Good and bad pants}
Essentially, a pants group is $(R, \epsilon)$-good (or bad) if there are $(R, \epsilon)$-well-connected (or badly-connected) tripods that generate the group. 

\begin{defn}
    We say a pair of pants $f\colon \Pi_0\to M$ is $(R,\epsilon)$-\emph{good}, if there are two $(R,\epsilon)$-well-connected framed tripods in $M$ with good connections $\gamma_i$, such that $\gamma_{i+1}\gamma_{i-1}^{-1}$ is freely homotopic to $f(C_i)$ for $i\in\Z/3\Z$. The $(R,\epsilon)$-\emph{bad pants} are defined the same way, but with bad connections. For convenience, we call a pants $(R,\epsilon)$-\emph{cuff-good}, if each cuff of this pants is an $(R,\epsilon)$-good curve.
\end{defn}

By Proposition 4.4 in \cite{Ham15}, we know that good pants have good boundary curves. Actually the same property holds for bad pants, which is not hard to verify, so we have the following proposition. 

\begin{prop}\label{cuffofgoodpants}
    There exists a universal constant $\chi>0$ such that any $(R,\epsilon)$-good or bad pants is $(R,\chi\epsilon)$-cuff-good.
\end{prop}

\begin{rem}
    Although Hamenst{\"a}dt ruled out the case of $\SO(2m,1)$ from the main theorem of \cite{Ham15}, those results before Section 7 in that paper work well for all $\SO(n,1)$ with $n\geq3$ as indicated in Section 7 of that paper.
\end{rem}

\subsubsection{The doubling trick}\label{sub:doubling}

This subsection follows Section 2.5 in \cite{KW21}. One can also refer to Section 4.7 in \cite{Rao25}.

Suppose now that we have a family of pants $\mathcal{F}=\{f_i\colon P_i\to M\}$, where $P_i$'s are orientable. We take two copies of each $f_i\colon P_i\to M$ with opposite orientation and denote them by $f_i^+\colon P_i^+\to M$ and $f_i^-\colon P_i^-\to M$. We then obtain a family of oriented pants $2\mathcal{F}$. 
(We think of the $P_i$ as distinct (and disjoint) pants, so $\bigcup \FF$ is a disjoint union of pants, and likewise for $\bigcup 2\FF$.)

For each geodesic $\gamma\in M$, let $\mathcal{F}(\gamma)$ be the subset of $\cup_i\partial P_i$ consisting of those elements that are mapped to the unoriented $\gamma$ in $M$. $\mathcal{F}(\gamma)$ could be empty if $\gamma$ is not the image of any boundary component. We define $2\mathcal{F}(\gamma)$ in a similar way with $2\mathcal{F}$, by taking two copies of each $\alpha\in\mathcal{F}(\gamma)$ with opposite orientation, which is a subset of $(\cup_i\partial P_i^+)\cup(\cup_i\partial P_i^-)$. For $\alpha\in\mathcal{F}(\gamma)$, there is $\alpha^+\in2\mathcal{F}(\gamma)$ that is mapped to $\gamma$ and $\alpha^-\in2\mathcal{F}(\gamma)$ that is mapped to $\gamma^{-1}$. Now suppose that there is a self-bijection $\sigma_\gamma$ on $\mathcal{F}(\gamma)$. We then can define an involution $\tau_\gamma$ on $2\mathcal{F}(\gamma)$ by
\begin{equation*}
    \tau_\gamma(\alpha^+)=(\sigma_\gamma(\alpha))^-
\end{equation*}
and
\begin{equation*}
    \tau_\gamma(\alpha^-)=(\sigma_\gamma^{-1}(\alpha))^+.
\end{equation*}
By combining each involution for all $\gamma$, we get an involution $\tau$ on $(\cup_i\partial P_i^+)\cup(\cup_i\partial P_i^-)$. Then we identify $\alpha$ and $\tau(\alpha)$ for each $\alpha\in(\cup_i\partial P_i^+)\cup(\cup_i\partial P_i^-)$, and we will get a closed surface $\mathcal{S}$ and a map $f\from\mathcal{S}\to M$ by joining $P^\pm_i$'s and $f^\pm_i$'s via $\tau$.

\subsection{Steiner graph of pants}
We first recall the result of Fermat point for hyperbolic triangles from Section 5 in \cite{TX18}. Suppose $\triangle ABC$ is a hyperbolic triangle with three angles $\angle A,\angle B,\angle C$. When $\angle A,\angle B,\angle C<2\pi/3$, then there is a unique Fermat point $P$ inside the triangle such that $|PA|+|PB|+|PC|$ is minimal. When one of three angles is no less than $2\pi/3$, then the vertex corresponding to this angle is the unique Fermat point. In the former case, 
\begin{equation}
    \angle APB=\angle BPC=\angle CPA=2\pi/3.
\end{equation}

For $A, B \in M$,
we let $\pi_1(M,A,B)$ denote the set of homotopy classes (rel endpoints) of paths in $M$ that connect $A$ and $B$. Then $\pi_1(M, A, A)$ can be identified with $\pi_1(M, A)$, and we can form $\alpha^{-1} \in \pi_1(M, B, A)$ and $\alpha \beta \in \pi_1(M, A, C)$ when $\alpha \in \pi_1(M, A, B)$ and $\beta \in \pi_1(M, B, C)$.

We can now be more precise with the definition of a pants $\Pi$ in $M$.
We can take a base point $*$ for our model pants $\Pi_0$ and draw three paths $\gamma_i$ from the base point to itself, such that each $\gamma_i$ is freely homotopic to a distinct cuff. Then the $\gamma_i$'s generate $\pi_1(\Pi_0)$, and $\gamma_0\gamma_1\gamma_2 = 1$. 
Given any $P \in M$, a pair of pants $\Pi$ in $M$ is determined by a representation $\rho\from\pi_1(\Pi_0, *) \to \pi_1(M, P)$;
we say $\rho$ and $\rho'\from\pi_1(\Pi_0, *) \to \pi_1(M, P')$ determine the same pants if there is $\eta \in \pi_1(M, P, P')$ such that $\rho'(\gamma_i) = \eta^{-1} \rho(\gamma_i) \eta$ for each $i$.

Let $G(\Pi, M)$ denote the space of triples $(P, Q, (\eta_i))$,
where $P, Q \in M$,
each $\eta_i$ is a geodesic segment from $P$ to $Q$,
and the representation from $\pi_1(\Pi_0,*)$ to $\pi_1(M, P)$
given by $\gamma_{i+1} \mapsto \eta_i\eta^{-1}_{i+1}$ determines $\Pi$.
We define a function $l\from G(\Pi, M) \to [0, \infty)$ by
\begin{equation}
    l(P,Q, (\eta_i))=|\eta_1|+|\eta_2|+|\eta_3|.
\end{equation}
Since $G(\Pi, M)$ is locally compact and $l$ is continuous and proper, the function $l$ has a minimum at some $(P,Q, (\eta_i))\in G(\Pi, M)$. 
\begin{thm}
    The function $l$ has a unique local minimum (which is then the global minimum). 
\end{thm}
\begin{proof}
    Let us fix a base point $(P_*, Q_*, (\eta_{i, *}))$ in $G(\Pi, M)$, and let $\tilde M$ be the universal cover of $M$. We can think of $\tilde{M}$ as the space of all paths starting at $P_*$ up to homotopy rel endpoints; we can do the same with $P_*$ replaced by $Q_*$. Therefore we have
    \begin{equation}
        \tilde{M}\times\tilde{M}\cong \{(P, \alpha,Q,\beta):P,Q\in M, \alpha\in\pi_1(M,P_*,P),\beta\in\pi_1(M,Q_*,Q)\},
    \end{equation}
     We then can define the following map
    \begin{equation}
        \begin{aligned}
            \iota\colon \tilde{M}\times\tilde{M}&\to G(\Pi, M)\\
            (P,\alpha,Q,\beta)&\mapsto (P,Q,([\alpha^{-1}\eta_{i,*}\beta])),
        \end{aligned}
    \end{equation}
    where $[\gamma]$ is the geodesic segment that is path-homotopic to $\gamma$. We claim that $\iota$ is a bijection.

    If $(P,\alpha_1,Q,\beta_1)$ and $(P,\alpha_2,Q,\beta_2)$ have the same image, then 
    \begin{equation}
        \alpha_1^{-1}\eta_{i,*}\beta_1=\alpha_2^{-1}\eta_{i,*}\beta_2\in\pi_1(M,P,Q),
    \end{equation}
    for $i\in\Z/3\Z$. Thus $\eta_{i,*}^{-1}\alpha_2\alpha_1^{-1}\eta_{i,*}\beta_1\beta_2^{-1}$ is null-homotopic. Applying this for $i$ and $i+1$ and cancelling $\beta_1\beta_2^{-1}$, we obtain
    \begin{equation}
        \eta_{i,*}^{-1}\alpha_2\alpha_1^{-1}\eta_{i,*}
        =\eta_{i+1,*}^{-1}\alpha_2\alpha_1^{-1}\eta_{i+1,*},
    \end{equation}
    as elements in $\pi_1(M,Q_*)$. We then know that $\alpha_2\alpha_1^{-1}$ commutes with $\eta_{i,*}\eta^{-1}_{i+1,*}$ up to homotopy, as elements in $\pi_1(M,P_*)$. 
    Therefore $\alpha_2\alpha_1^{-1}$ lies in the centralizer of a pants group that determines $\Pi$, and is therefore trivial. Hence 
    $\alpha_1=\alpha_2$. 
    Similarly, $\beta_1=\beta_2$. Hence $\iota$ is injective.

    On the other hand, for any $(P,Q,(\eta_i))\in G(\Pi,M)$, the representation $\gamma_i\mapsto\eta_i\eta_{i+1}^{-1}$ determines the same pants $\Pi$ as the representation $\gamma_i\mapsto\eta_{i,*}\eta_{i+1,*}^{-1}$ which is given by $(P_*,Q_*,(\eta_{i,*}))$. Thus there exists $\alpha\in\pi_1(M,P_*,P)$ such that
    \begin{equation}\label{alphahomotopy}
        \eta_i\eta^{-1}_{i+1}=\alpha^{-1}\eta_{i,*}\eta^{-1}_{i+1,*}\alpha,
    \end{equation}
    for $i\in\Z/3\Z$. Let $\beta=\eta_{i,*}^{-1}\alpha\eta_i\in\pi_1(M,Q_*,Q)$. Then by \eqref{alphahomotopy}, we know $\beta$ does not depend on the choice of $i$. Moreover, we have $\eta_i=\alpha^{-1}\eta_{i,*}\beta$ for each $i\in\Z/3\Z$, which implies $\iota(P,\alpha,Q,\beta)=(P,Q,(\eta_i))$. Hence $\iota$ is also surjective, and then the claim is proved.
    
    

     We then can pull $l$ back to $\tilde M \times \tilde M$, and the unique local minimum is implied from the observation that $l$ is strictly convex (with respect to the product metric): this is because $l$ is the sum of three distances, and distance function is strictly convex on $\tilde M \times \tilde M$ (except on the diagonal, where it is identically 0).
\end{proof}
We observe that one of the $\eta_i$ may have length 0 if $P$ and $Q$ are the same point, but this can hold for at most one $\eta_i$. The minimum length graph with vertices $P,Q$ and connections $\eta_i$ is called the \emph{Steiner graph} of $\Pi$, and it is called \emph{degenerate} if one of the $\eta_i$ has length 0. We then prove that the Steiner graph of cuff-good pants is not degenerate and yields two tripods.

\begin{prop}\label{SteinerGraph}
    There exists $R_0>0$ such that for all $R>R_0$ and $0\le\epsilon<1$, the following holds: Suppose $\Pi$ is an $(R,\epsilon)$-cuff-good pants in a closed hyperbolic $n$-manifold $M$. Then the Steiner graph of $\Pi$ is not degenerate. Moreover, the angle between any two connections is $2\pi/3$ at each vertex of the Steiner graph.
\end{prop}

\begin{proof}
    Let $P,Q$ be the two vertices and $\eta_i$, $i\in\Z/3\Z$ be the edges of the Steiner graph of $\Pi$. We let $\tilde{Q}\in\H^n$ be a lift of $Q$ and $\tilde\eta_i$ be the lift of $\eta_i$ such that $\tilde\eta_i$ connects $\tilde{Q}$ and $\tilde{P}_i$, where $\tilde{P}_i$ are three different lifts of $P$. Since $l$ reaches its minimal value at $(P,Q,(\eta_i))$, then by the above result, $\tilde{Q}$ is the Fermat point of $\triangle \tilde{P}_1\tilde{P}_2\tilde{P}_3$. 
    
    If $(P,Q,(\eta_i))$ is degenerate, then $\tilde{Q}$ coincides with some $\tilde{P}_i$, and 
    \begin{equation}
        \begin{aligned}
            l(P,Q,(\eta_i)) &=\|\tilde{P}_i\tilde{P}_{i+1}\|+\|\tilde{P}_i\tilde{P}_{i+2}\|
            \ge l(\gamma_{i+2})+l(\gamma_{i+1})\\
            &\ge 2(2R-2\epsilon)=4R-4\epsilon.
        \end{aligned}
    \end{equation}
    Here, $\|XY\|$ denotes the hyperbolic distance between $X$ and $Y$ for $X,Y\in\H^n$.

    However, the function $l$ can actually attain smaller values. Let $\alpha_i$ be the short orthogeodesic between $\gamma_{i+1}$ and $\gamma_{i+2}$,
    then by hyperbolic trigonometry, there exists universal constant $C$ such that 
    \begin{equation}
        C^{-1}e^{-R/2}<l(\alpha_i)<Ce^{-R/2}.
    \end{equation}
    Suppose $P'=\alpha_0\cap\gamma_1$, and we let $\tilde{P}'_1, \tilde{P}'_2, \tilde{P}'_3$ be three lifts of $P'$ such that the projection of $\tilde{P}'_i\tilde{P}'_{i+1}$ is freely homotopic to $\gamma_{i+2}$. We then let $\tilde{Q}'$ be the Fermat point of the hyperbolic triangle $\triangle\tilde{P}'_1\tilde{P}'_2\tilde{P}'_3$, and $Q'$ be the projection of $\tilde{Q}'$. Since $P'\in\gamma_1$, then the projection of $\tilde{P}'_2\tilde{P}'_{0}$ in $M$ is $\gamma_1$, so $\|\tilde{P}'_2\tilde{P}'_{0}\|=l(\gamma_1)$. By hyperbolic geometry, we know that there exists a universal constant $C>0$ such that
    \begin{equation}\label{P_iP_j}
        \begin{aligned}
            2R-2\epsilon\le l(\gamma_2)&\le\|\tilde{P}'_0\tilde{P}'_{1}\|\le2R+C,\\
            2R-2\epsilon\le l(\gamma_0)&\le\|\tilde{P}'_1\tilde{P}'_{2}\|\le3R+C.
        \end{aligned}
    \end{equation}
    \begin{figure}[htpb]
        \centering
        \includegraphics[width=0.8\textwidth]{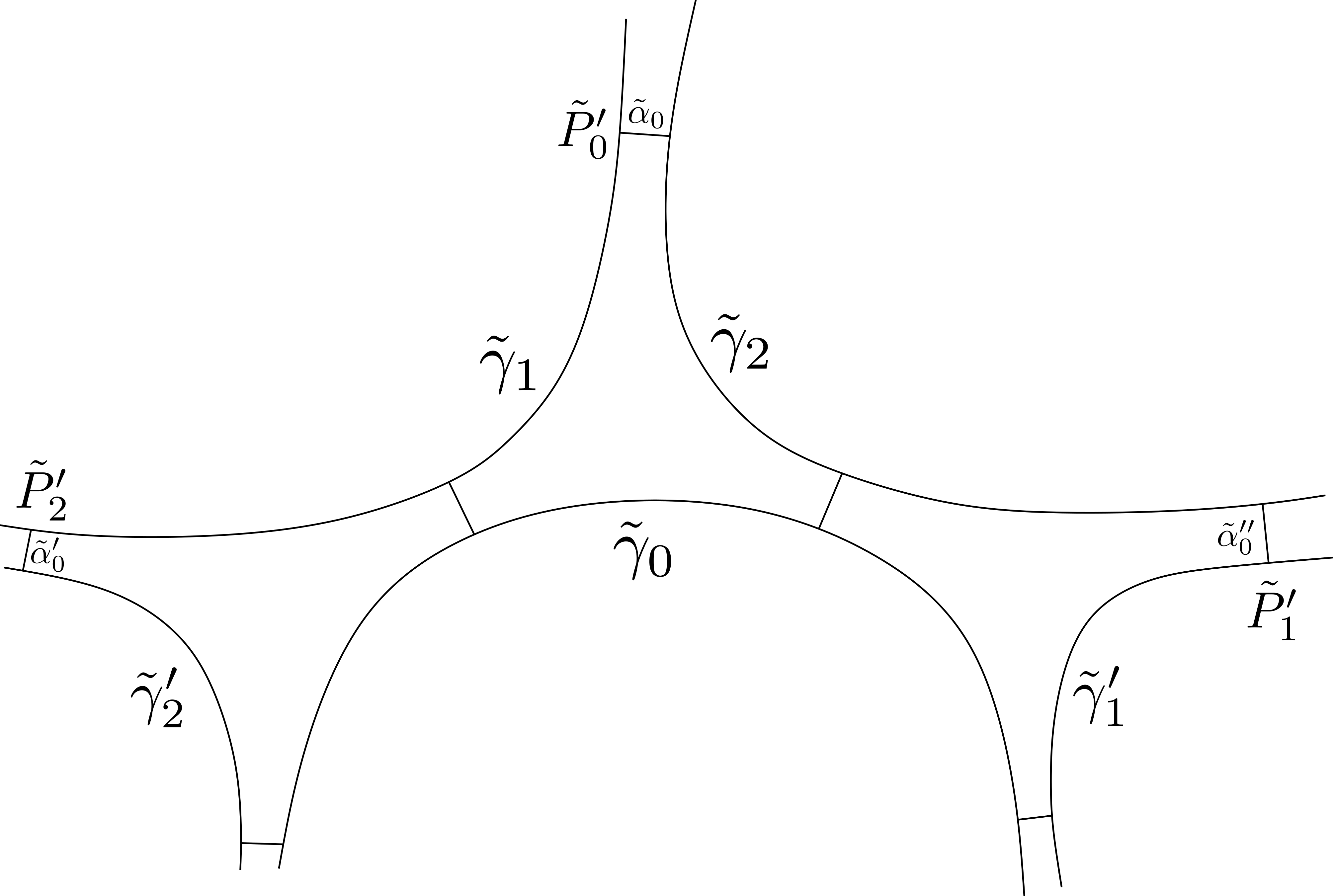}
        \caption{Lifts of $P'$ in the universal cover}
        \label{fig:Steiner}
    \end{figure}
    Therefore when $R$ is large enough, each angle of $\triangle\tilde{P}'_1\tilde{P}'_2\tilde{P}'_3$ is less than $2\pi/3$. Hence, the Fermat point $\tilde{Q}'$ is in the interior. Furthermore, there is a universal constant $C'$ such that 
    \begin{equation}\label{Fermatpointvalue}
        \|\tilde{Q}'\tilde{P}'_0\|+\|\tilde{Q}'\tilde{P}'_1\|+\|\tilde{Q}'\tilde{P}'_2\|
        \le
        \frac{1}{2}\left(\|\tilde{P}'_0\tilde{P}'_1\|+\|\tilde{P}'_1\tilde{P}'_2\|+\|\tilde{P}'_2\tilde{P}'_0\|\right)
        +C'
    \end{equation}
    Thus when $R$ is sufficiently large, by \eqref{P_iP_j} and \eqref{Fermatpointvalue}, we have
    \begin{equation}
        \begin{aligned}
            l(P',Q,(\eta'_i))
            &=\|\tilde{Q}'\tilde{P}'_0\|+\|\tilde{Q}'\tilde{P}'_1\|+\|\tilde{Q}'\tilde{P}'_2\|\\
            &\le\frac{1}{2}\left((2R+2\epsilon)+(2R+C)+(3R+C)\right)+C'\\
            &=3.5R+\epsilon+C+C'<4R-4\epsilon,
        \end{aligned}
    \end{equation}
    where $\eta'_i$ is the projection of $\tilde{P}'_{i+1}\tilde{P}'_{i+2}$. That is contradictory to our assumption that $(P,Q,(\eta_i))$ was the Steiner graph. Therefore the Steiner graph of $\Pi$ is not degenerate.

    Furthermore, we know that $\tilde{Q}$, which is the Fermat point of $\triangle \tilde{P}_1\tilde{P}_2\tilde{P}_3$, is in the interior of $\triangle \tilde{P}_1\tilde{P}_2\tilde{P}_3$. Thus $\angle\tilde{P}_i\tilde{Q}\tilde{P}_{i+1}=2\pi/3$, for $i\in\Z/3\Z$. Thus the angle between any two connections at $Q$ is $2\pi/3$. We then know that the result holds at $P$ as well by the symmetry of $P$ and $Q$.
\end{proof}

\begin{rem}
    There are actually three (reduced) graphs with the homotopy type of a pants: the ``theta", which is what we are considering, the ``barbell'', where there are segments from each point $P, Q$ to itself, and a segment connecting them, and the ``figure eight'', which can be thought of as either a degenerate theta or a degenerate barbell, and which therefore forms the boundary between them. One can then prove that there is a unique local minimum on the total space of graphs. 
\end{rem}

\subsection{Good cuffs implies good or bad pants}
In this section, we prove a surprising fact, that the converse of Proposition \ref{cuffofgoodpants} is also true, up to an adjustment of constant.

For $A,B\in \SO(n-1)$, we say $A$ is $\epsilon$-close to $B$ if $AB^{-1}$ is $\epsilon$-close to $\mathbf{e}$. Let $Z=\mathbf{diag}(-1,1,\ldots,1)\in \mathrm{O}(n-1)$. The conjugacy by $Z$ can be restricted to a self-map of $\SO(n-1)$, and we denote it by $\phi\colon \SO(n-1)\to \SO(n-1)$. Let
\begin{equation*}
    \mathcal{Q}=\left\{\begin{pmatrix}
        1 & 0 \\
        0 & X
    \end{pmatrix}\colon X\in SO(n-2)\right\}
\end{equation*}
be a subgroup of $\SO(n-1)$, which is isomorphic to $\SO(n-2)$. Then the elements in $\mathcal{Q}$ are fixed by $\phi$.

\begin{lem}\label{2.7}
    Suppose $A\in \SO(n-1)$. Then there exists $Q\in \mathcal{Q}$, 
    and $U \in \{ \mathbf{e}, \mathbf{diag}(-1,-1,1,\ldots,1) \}$ such that 
    $$d(AQ, U)  \le  \frac12 d(A, \phi(A)).$$
\end{lem}

\begin{proof}
Let $e_1 = (1, 0, \ldots, 0) \in \R^{n-1}$, so $Q \in \mathcal{Q}$ if and only if $Q e_1 = e_1$. 
For any $v \in S^{n - 2}$, 
we can take $z \in \{e_1, -e_1\}$ such that
$d(-v, Z v) = 2d(v, z)$.
We then have 
$$
d(A, \phi(A)) = d(A, ZAZ) = d(AZ, ZA) 
\ge d(-A e_1, Z A e_1) = 2 d(A e_1, U e_1)$$
for either $U = \mathbf{e}$ or $U = \mathbf{diag}(-1,-1,1,\ldots,1)$.
For convenience of notation we note that $U^{-1} = U$. 

We can then find $R \in \SO(n -1)$ such that $R A e_1 = U e_1$,
 and $d(R, \mathbf{e}) =  d(A e_1, U e_1)$.
Then $URA e_1 = e_1$, so we can find $Q \in \mathcal{Q}$ such that $URA = Q^{-1}$. 
Then $AQ = R^{-1}U$, so $d(AQ, U) = d(R, \mathbf{e}) \le \frac12 d(A, \phi(A))$. 
\end{proof}

Let
\begin{equation*}
    R(\theta)=\begin{pmatrix}
        \cos\theta & \sin\theta & & &\\
        -\sin\theta & \cos\theta & & &\\
        & & 1 & &\\
        & & & \ddots &\\
        & & & & 1
    \end{pmatrix}\in \SO(n), \mathrm{for}\ \theta\in[0,2\pi];
\end{equation*}
then the right action of $R(\theta)$ on a n-frame is the rotation by $\theta$ in the plane spanned by the first two vectors, and $Z$ can be viewed as the restriction of $R(\pi)$ on the space spanned by the last $n-1$ vectors of the frame. For any $t\in\R$, let $G(t)$ be the right action of the frame flow along the first vector of the frame by distance $t$. We notice that $R(\pi)^2=\mathrm{Id}$, and $R(\pi)$ conjugates $G(t)$ to $G(-t)$, for any $t$.

The following lemma generalizes the inefficiency theory of hyperbolic geometry to dimension $n$. We will skip the proof and one can refer to \cite{KM15}, \cite{LM15}, or \cite{Rao25} for more details.

\begin{lem}[Inefficiency of Broken geodesic segments]\label{ineff}
    For any $0<\epsilon<\pi$, there exist $C,L,m>0$ such that for any $\theta\in(\epsilon-\pi,\pi-\epsilon)$: Suppose $t_1,t_2>L$, then there exist $Y_1,Y_2\in \SO(n)$ that are $Ce^{-L}$-close to $\mathbf{e}$, and $t_1+t_2>t>t_1+t_2-m$ such that
    \begin{equation*}
        G(t_1)R(\theta)G(t_2)=Y_1G(t)Y_2.
    \end{equation*}
\end{lem}

The next theorem is the reverse statement of Proposition \ref{cuffofgoodpants}. It shows that the length and monodromy of cuffs can provide enough information to tell whether the pants is good or not. We can then consider the set of all $(R, \epsilon)$-cuff-good (and good) pants and use our estimates in Section \ref{sec:counting} to prove in Section \ref{sec:final-proof} that we can match these pants to form an incompressible surface. 

\begin{thm}\label{2.9}
There is a universal constant $C$ such that if $\Pi$ is a pants group that is $(R, \epsilon)$-cuff-good, then $\Pi$ is either $(R, C\epsilon)$-good or $(R, C\epsilon)$-bad. 
\end{thm}

\begin{proof}
    Suppose the three cuffs of an oriented pants $\Pi$ in $M$ are $\gamma_0,\gamma_1$ and $\gamma_2$ that are $(R,\epsilon)$-good. We proceed by the following steps:
    \begin{enumerate}
        \item 

        \textbf{Use the Steiner graph of $\Pi$ to get two (triply) connected tripods.}

        Let $G$ be the Steiner graph of $\Pi$ with vertices $p,q$ and edges $\eta_0,\eta_1,\eta_2$ connecting $p$ and $q$, such that $\eta_{i+1} \eta_{i-1}^{-1}$ is freely homotopic to $\gamma_{i}$, for $i\in\Z/3\Z$, and we denote by $l_i$ the length of $\eta_i$. Without loss of generality, we assume for each $i\in\Z/3\Z$, $\eta_i$ is parameterized by arc length and oriented from $p$ to $q$. By a mild abuse of notation, when $z = \eta_i(t)$, we will write $\eta'(z)$ instead of $\eta'(t)$. This is unambiguous because each $\eta_i$ is injective.
        By Prop \ref{SteinerGraph}, the angle between $\eta'_i(p)$ and $\eta'_j(p)$ is $2\pi/3$ if $i\neq j\in\Z/3\Z$, and the same is also true for $-\eta'_i(q)$ and $-\eta_j'(q)$.
        Thus $P=(p,\eta'_0(p),\eta'_1(p),\eta'_2(p))$ and $Q=(q,-\eta'_0(q),-\eta'_{-1}(q),-\eta'_{-2})(q)$ are two tripods generating $\Pi$. Then we only need to find frames at $p$ and $q$ such that the monodromy along each $\eta_i$ is $C\epsilon$-close to $\mathbf{e}$ or $\mathbf{diag}(1,-1,-1,1,\ldots,1)$, for some constant $C$.
        
        \item 
        \textbf{Choose frames for the orthogonal plane at each tripod, to get complete frames for $\R^{n}$.}
        
        We first choose arbitrary positive frames $E$ and $F$ at $p$ and $q$ in the orthogonal complement of the 2-planes $V_p$ and $V_q$ determined by the tripod $P$ and $Q$, respectively. Then $P_0=(P,E)$ and $Q_0=(Q,F)$ are two framed tripods that generate $\Pi$. For each $i\in\Z/3\Z$, there are vectors $u_i$ at $p$ and $v_i$ at $q$ such that $(\eta'_i(p),u_i)$ and $(-\eta'_i(q),v_i)$ are positive frames of $V_p$ and $V_q$, respectively. Let $E_i=(\eta'_i(p),u_i,E)$ and $F_i=(-\eta'_i(q),v_i,F)$ be the frames at $p$ and $q$. Then the $E_i$ are related by rotations on $V_p$, and the case is similar for the $F_i$, where $R(2\pi/3)$ sends $E_i$ and $F_i$ to $E_{i+1}$ and $F_{i-1}$, respectively.
        
        \item 
        \textbf{For each connection, push these two frames to the midpoint of the connection and obtain a transformation $X_i$ between the two new frames. Then write the monodromy of each cuff as a word in these $X_i$.}

        Let $x_i$ be the midpoint of $\eta_i$, for $i\in\Z/3\Z$. We then push forward $E_i$ and $F_i$ to $x_i$ by parallel transport, and the frames we get at $x_i$ are denoted by $E_i(x_i)$ and $F_i(x_i)$. Then there is a unique $X_i\in \SO(n)$ such that the right action of $R(\pi)X_i$ sends $E_i(x_i)$ to $F_i(x_i)$. We notice that the right action of $R(\pi)$ on $E_i(x_i)$ will send $\eta'_i(x_i)$ to $-\eta'_i(x_i)$, therefore $X_i$ preserves $-\eta'_i(x_i)$. 

        Now we consider the conjugacy class of the element in $\mathcal{G}$ of the broken geodesic $\eta_{i+1}\eta^{-1}_{i+2}$, by applying a sequence of right actions to the frame $E_{i+1}(x_{i+1})$ and getting back to itself:
        \begin{equation}
        \begin{aligned}
            E_{i+1}(x_{i+1}) &\xrightarrow{R(\pi)X_{i+1}} F_{i+1}(x_{i+1}) \xrightarrow{G(-l_{i+1}/2)} F_{i+1} \xrightarrow{R(-2\pi/3)} F_{i+2} \xrightarrow{G(l_{i+2}/2)} F_{i+2}(x_{i+2})\\
            &\xrightarrow{(R(\pi)X_{i+2})^{-1}} E_{i+2}(x_{i+2}) \xrightarrow{G(-l_{i+2}/2)} E_{i+2} \xrightarrow{R(-2\pi/3)} E_{i+1} \xrightarrow{G(l_{i+1}/2)} E_{i+1}(x_{i+1}).
        \end{aligned}
        \end{equation}
        Thus the element in $\mathcal{G}$ of the homotopy class of $\eta_{i+1}\eta^{-1}_{i+2}$ is, up to conjugacy, 
        \begin{equation*}
        \begin{split}
            &[R(\pi)X_{i+1} G\left(-\frac{l_{i+1}}{2}\right) R\left(\frac{-2\pi}{3}\right)  G\left(\frac{l_{i+2}}{2}\right) X_{i+2}^{-1} R(\pi)  \\
            &G\left(-\frac{l_{i+2}}{2}\right) R\left(\frac{-2\pi}{3}\right)  G\left(\frac{l_{i+1}}{2}\right)]\\
            =\ &[R(\pi)X_{i+1} R(\pi)G\left(\frac{l_{i+1}}{2}\right)R(\pi) R\left(\frac{-2\pi}{3}\right)  \\
            &G\left(\frac{l_{i+2}}{2}\right) X_{i+2}^{-1}  G\left(\frac{l_{i+2}}{2}\right) R(\pi) R\left(\frac{-2\pi}{3}\right)  G\left(\frac{l_{i+1}}{2}\right)]\\
            =\ &[R(\pi)X_{i+1} R(\pi)G\left(\frac{l_{i+1}}{2}\right) R\left(\frac{\pi}{3}\right)  \\
            &G\left(\frac{l_{i+2}}{2}\right) X_{i+2}^{-1}  G\left(\frac{l_{i+2}}{2}\right) R\left(\frac{\pi}{3}\right)  G\left(\frac{l_{i+1}}{2}\right)]
        \end{split}
        \end{equation*}
        When $R$ is large enough, by Lemma \ref{ineff}, we can find $Y_i^{\pm}\in\SO(n)$, which are $Ce^{-R/2}$-close to $\mathbf{e}$, and $l_i^\pm\in\R_+$, $i\in\Z/3\Z$ such that
        \begin{equation*}
            \begin{split}
                G\left(\frac{l_{i+1}}{2}\right)R\left(\frac{\pi}{3}\right)G\left(\frac{l_{i+2}}{2}\right)&=Y_{i+1}^+G(l_i^+)Y_{i+2}^+,\\
                G\left(\frac{l_{i+2}}{2}\right)R\left(\frac{\pi}{3}\right)G\left(\frac{l_{i+1}}{2}\right)&=Y_{i+2}^-G(l_i^-)Y_{i+1}^-.
            \end{split}
        \end{equation*}
        Then the above conjugacy class can be written as
        \begin{equation*}
            \begin{split}
                &[R(\pi)X_{i+1} R(\pi) Y_{i+1}^+G(l_i^+)Y_{i+2}^+ X_{i+2}^{-1}  Y_{i+2}^-G(l_i^-)Y_{i+1}^-]\\
                =\ &[G(l_i^+)Y_{i+2}^+ X_{i+2}^{-1}  Y_{i+2}^-G(l_i^-)Y_{i+1}^-R(\pi)X_{i+1} R(\pi) Y_{i+1}^+].
            \end{split}
        \end{equation*}
        By the definition of $X_i$ and $R(\pi)$, we know that $X_i$ and $R(\pi)X_{i} R(\pi)$ both preserve the first vector of  frames, and hence lie in $\mathcal{M}$.  We are now in a position to apply Lemma \ref{8-word}, where the $Y$'s take on the role of the $u$'s and $v$'s in that Lemma, and $X^{-1}_{i+2}$ and $R(\pi)X_{i+1}R(\pi)$ become  $m_1$ and $m_2$. As in the preamble to Lemma \ref{8-word}, 
        we conclude that the length of $\gamma_i$ can be approximated by $l_i^++l_i^-$, and the monodromy of $\gamma_i$ (starting from $x_{i+1}$) can be estimated by $X_{i+2}^{-1}R(\pi)X_{i+1} R(\pi)$, where the error terms are bounded by $e^{-R/2}$ times a universal constant.  

        \item
        \textbf{These transformations are right actions acting on the frames, but we can write the right actions as matrices and apply Lemma \ref{2.7}}.
        
        Let us now consider the element $X_{i+2}^{-1}R(\pi)X_{i+1} R(\pi)$ of $\SO(n) \subset \mathcal{G}$. From Step 3 we know that, for large $R$, it is $\epsilon$-close to the monodromy of $\gamma_i$ (which by our assumption, is $\epsilon$-close to $\mathbf{e}$),  and is therefore $2 \epsilon$-close to $\mathbf{e}$ . On the other hand we know that this right action preserves the first vector of the frame. As the $X_i$ also preserve the first vector, we can think of them as elements of $\SO(n-1)$; recalling that $Z$ is the restriction of $R(\pi)$ to the last $n-1$ vectors, we can then write $X_{i+2}^{-1}R(\pi)X_{i+1} R(\pi)$, as an element of $\SO(n-1)$, as 
        $X_{i+2}^{-1} \phi(X_{i+1})$. 
        Hence $X_{i+2}$ is $2\epsilon$-close to $\phi(X_{i+1})$, $i\in\Z/3\Z$, and we then know $X_i$ is $6\epsilon$-close to $\phi(X_i)$, and $X_i$ is $4\epsilon$-close to $X_j$, for $i\in\Z/3\Z$ and $j\neq i$. By Lemma \ref{2.7}, we can find $Q\in\mathcal{Q}$ such that $X_0Q$ is $3\epsilon$-close to either $\mathbf{e}$ or $\mathbf{diag}(-1,-1,1,\ldots,1)$. In the former case, we have $X_iQ$ is $7\epsilon$-close to $\mathbf{e}$, for all $i\in\Z/3\Z$. Since the right action of $Q$ commutes with the frame flow and the rotations $R(\theta)$, we can replace $X_i$ with $X_iQ$ through replacing $F$ with $F\cdot Q$. We then have a pair of well-connected framed tripods. 
        
        For the latter case, each $X_iQ$ is $7\epsilon$-close to $\mathbf{diag}(-1,-1,1,\ldots,1)$ as elements in $\SO(n-1)$. Since they preserve the first vector of the frame when considered as elements in $\SO(n)$, each $X_iQ$ is $7\epsilon$-close to $\mathbf{diag}(1,-1,-1,1,\ldots,1)$ in $\SO(n)$. Hence we obtain a pair of badly-connected framed tripods by the same reasoning. \qedhere
    \end{enumerate}
\end{proof}


\section{Quasi-uniformity for the measure on feet of good pants} \label{sec:counting}

\textbf{We will assume throughout this section} that $\gamma_0$ is a given $(R, \epsilon)$-good curve.
We want to find all third connections $\eta$ such that the pants generated by each $\eta$ is an $(R,\epsilon)$-cuff-good pants, but also a \emph{good} pants (i.e., not a \emph{bad} one). More specifically, we want to know the distribution of the (long) feet of these third connections, and roughly what measure we get on the average of a pair of feet. 



Given two unit normal vectors on $\gamma_0$ as two feet of a predicted third connection, we can estimate the volume of the space of monodromies for this third connection for which the two other cuffs will be good. This volume is zero unless the two unit normal vectors are close (when compared in either direction) and basically just depends on how close they are. We will make this much more explicit in this section.



In this way we obtain an effective estimate of the measure on the feet of the good third connections (first on the long feet, and then on their midpoint). With this, we will prove in Section \ref{sec:final-proof} that,
for each good curve $\gamma_0$, we can match the pants around $\gamma_0$ such that the feet for matching pants satisfy the usual condition of being on opposite sides of $\gamma_0$ (nearly antipodal points of the sphere) along with a shear by 1. We can then apply Hamenst{\"a}dt's sufficient condition to conclude that the resulting surface is $\pi_1$-injective. 

In Section \ref{sec:third-conn-length},
we estimate the lengths of the two new cuffs in terms of the distances (along $\gamma_0$) between the footpoints of the third connection. 
In Section \ref{sec:computingmonodromy}, 
we compute the monodromy of the two new cuffs in terms of right actions for frames, and then estimate these monodromies in terms of the distance between the feet of the third connections. In particular we see that the feet must be $3\epsilon/2$-close (after parallel transport) if the monodromies are good, and the mondromies are good if the feet are $\epsilon/4$-close. 
In Section \ref{sec:ave-feet},
we consider the \emph{average feet} of the long feet, and observe that these are close to the short feet.
Finally, 
in Section \ref{sec:volume},
we first compute the volume of the space of good potential third connections for a given pair of feet for the third connection, and then compute the volume of the space of good potential third connections for a given average foot. We then conclude, in Proposition \ref{original-counting},
that the measure $\mu_a$ for average feet for good potential third connections is bounded above and below in terms of Lebesgue measure. 

This is the central estimate, which we will apply in Section \ref{sec:final-proof} as follows. By general principles the measure for average feet for good potential third connections is invariant by the ``twist by 1 and antipodal map'' that we use for matching. This is essentially the \emph{estimated} foot measure, but it's not exactly the \emph{actual} foot measure, which we need to be invariant by the matching map, \emph{up to a small error in the sense of mass transport}.  
But this form of invariance follows from the quasi-uniformity of the estimated foot measure, along with the invariance by the matching map. So we can perform the matching and then conclude that the resulting surface is $\pi_1$-injective.




\subsection{Lengths of cuffs determined by a third connection}\label{sec:third-conn-length}
In this subsection we describe the region of pairs of endpoints of a putative third connection that generates cuffs with good lengths. We then estimate the area of this region. 

We start with another element of the Theory of Inefficiency, which is like Lemma \ref{ineff}, except for a closed broken geodesic, instead of a broken geodesic segment. This lemma helps us estimate the length of a third connection in terms of the length of cuffs in good pants. The same estimate in the two and three dimensional cases were used in \cite[Section 3.3]{KM15} and \cite[Section 3.3]{KW21}, respectively, which is actually a special case of Lemma \ref{ineff}. This lemma helps us estimate the length of a third connection in terms of the length of cuffs in good pants. The same estimate in the two and three dimensional cases were used in \cite[Section 3.3]{KM15} and \cite[Section 3.3]{KW21}, respectively.

\begin{lem}\label{lem:third-conn-length}
    There exists constant $C,L_0>0$ such that for all $L>L_0$, the following holds: Suppose $\gamma^+,\gamma^-$ are two oriented geodesic segments in $M$ which meet orthogonally at both endpoints satisfying that $\gamma^+\gamma^-$ is an oriented broken geodesic. If $l(\gamma^i)>L$, $i=+,-$, then
    \begin{equation*}
        |l(\gamma)-l(\gamma^+)-l(\gamma^-)+2\ln2|<Ce^{-L},
    \end{equation*}
    where $\gamma$ is the closed geodesic homotopic to $\gamma^+\gamma^-$.
\end{lem}

\begin{rem}
    The monodromy of geodesics in hyperbolic 3-manifolds can be regarded as the imaginary part of the complex length, which cannot be defined properly in higher-dimensional cases. Therefore we need to study the length and the monodromy of third connections separately. 
\end{rem}

Let $\eta$ be a third connection for a cuff $\gamma_0$ of an $(R,\epsilon)$-good pants, and suppose $\eta$ subdivides $\gamma_0$ into two arcs $\alpha_1$ and $\alpha_2$. We let $u_i$ be the length of $\alpha_i$ for $i=1,2$, and $w$ be the length of $\eta$. Then $u_1+u_2=l(\gamma_0)$. Let $L_{R,\epsilon}(\gamma_0)$ be the set of all pairs $(x,y)$ such that
\begin{equation}
    |x+y-2\ln2-2R|<2\epsilon
\end{equation}
and
\begin{equation}
    |(l(\gamma_0)-x)+y-2\ln2-2R|<2\epsilon.
\end{equation}
Then the square region $L_{R,\epsilon}(\gamma_0)\subset\R^2$ can be used to approximate the set of $(u_1,w)$ for all third connections of good pants, up to an exponentially small error term. Suppose that the measure on $\R\times\R$ is given by $e^{2y}dxdy$, then the area of $L_{R,\epsilon}(\gamma_0)$ is
\begin{equation}\label{eq:diamondregion1}
    \Vol(L_{R,\epsilon}(\gamma_0))=8e^{4R-l(\gamma_0)}(e^{2\epsilon}-e^{-2\epsilon})^2\sim 128\epsilon^2 e^{4R-l(\gamma_0)},
\end{equation}
as $\epsilon\to0$.

\subsection{Monodromy estimate and feet distance}\label{sec:computingmonodromy}
In this subsection, we first compute the monodromies of the other two cuffs $\gamma_1, \gamma_2$ from the frames at the footpoints of $\eta$, and then determine the ``good region'' of monodromies for the third connections that generate cuff-good pants. We assume in this subsection that the endpoints of the third connection determine good \emph{lengths} for the new cuffs, as described in Section \ref{sec:third-conn-length}, so we can concentrate on the monodromies. 


Let $\eta$ be a third connection for $\gamma_0$, where the endpoints of $\eta$ are roughly $R$ apart (along $\gamma_0$) and the length of $\eta$ is about $R + 2 \ln 2$. We want to determine when the remaining two cuffs of the pants determined by $\gamma_0$ and $\eta$ are $(R,\epsilon)$-good curves. Let $A,B$ be the points where $\gamma_0$ meets $\eta$ such that $\eta$ is oriented from $A$ to $B$, and let $S_i$ be the unit normal sphere of $\gamma_0$ at $i$ for $i=A,B$. Let $\gamma_0$ be divided by $\eta$ into $\alpha_1,\alpha_2$ such that $\gamma_0=\alpha_1\alpha_2^{-1}$, where $\alpha_1$ and $\alpha_2$ are both oriented from $A$ to $B$. Let $\gamma_1$ be the oriented cuff homotopic to $\eta\alpha_1^{-1}$ and $\gamma_2$ be the oriented cuff homotopic to $\eta^{-1}\alpha_2$. Let $y$ be the midpoint of $\eta$, and $x_i$ be the midpoint of $\alpha_i$ for $i=1,2$.

We now want to set up a system of frames that will allow us to write the monodromies for the $\gamma_i$ as a product of right actions. 
At point $A$, the unit tangent vectors $\gamma'_0(A)$ of $\gamma_0$ and $\eta'(A)$ of $\eta$ form a 2-plane $P_A$ in the tangent space $T_A M$, 
and we temporarily choose an arbitrary frame $E$ of the orthogonal complementary subspace of $P_A$ in $T_A M$. Let $E_0=(\eta'(A),-\gamma'_0(A),E)$, $E_1=(\gamma'_0(A),\eta'(A),E)$ and $E_2=(-\gamma'_0(A),-\eta'(A),E)$ be three frames at $A$. 
Similarly, we choose an arbitrary frame $F$ of the orthogonal complementary subspace of $P_B$ in $T_B M$, where $P_B$ is the 2-plane spanned by $\gamma'_0(B)$ and $\eta'(B)$. Then we have three frames $F_0=(-\eta'(B),-\gamma'_0(B),F)$, $F_1=(-\gamma'_0(B),\eta'(B),F)$ and $F_2=(\gamma'_0(B),-\eta'(B),F)$ at $B$.

We apply the frame flow to $E_0$ and $F_0$ until they arrive at $y$, and denote the resulting frames by $E_0(y)$ and $F_0(y)$. There then exists a unique $Y\in\SO(n)$ such that the right action by $R(\pi)Y$ sends $E_0(y)$ to $F_0(y)$. 
Similarly, we apply the frame flow to $E_i,F_i$, for $i=1,2$, until they arrive at $x_i$; 
we again have $X_i\in\SO(n)$ such that $E_i(x_i) \cdot R(\pi)X_i = F_i(x_i)$. 
From $Y$ and the $X_i$, we can compute the conjugacy class of the element in $\mathcal{G}$ for each $\gamma_i$, $i\in\Z/3\Z$;  it can be expressed as a composition of right actions as in the proof of Theorem \ref{2.9}. Then as in that proof, 
\begin{enumerate}
\item[(*)]
the monodromies for $\gamma_0$,$\gamma_1$ and $\gamma_2$
are given,
up to conjugacy and an $e^{-R/2}$ error term,
by $X_2^{-1}\phi(X_1)$, $X_1^{-1}\phi(Y)$ and $Y^{-1}\phi(X_2)$. 
\end{enumerate}
Moreover $Y^{-1}\phi(X_2)$ is conjugate to $(X_2^{-1}\phi(Y))^{-1}$; from here on we will consider the $X_i$ and $Y$ as elements in $\SO(n-1)$ since they all preserve the first vector in the frame.
This means (ignoring for the moment the very small error term),
the resulting pants is cuff-good if and only if  $X_1^{-1}\phi(Y)$ and  $X_2^{-1}\phi(Y)$ are $\epsilon$-close to $\mathbf{e}$. 

Therefore the space of
good third connections, with given good feet, is equal to the space of $Y$ for which $\phi(Y)$ is $\epsilon$-close to $X_1$ and $X_2$.
Since $\phi$ is volume preserving,  our goal is to \textbf{estimate the volume of the intersection of $\epsilon$-neighborhoods of $X_1$ and $X_2$ (essentially determined by the distance between $X_1$ and $X_2$ in $\SO(n-1)$) whenever $X_2^{-1}\phi(X_1)$ is $\epsilon$-close to $\mathbf{e}$}. 


Let $\Lambda$ be the monodromy of $\gamma_0$. We define the spherical distance (or the angle) between $\eta'(A)$ and $-\eta'(B)$ along $\alpha_1$, denoted by $d^1_{S^{n-2}}(\eta'(A),-\eta'(B))$, as the angle between these two unit normal vectors when we bring them to $x_1$ by parallel transport. If we bring them to $x_2$, the distance will equal the distance between $\Lambda^{-1}(\eta'(A))$ and $-\eta'(B)$ along $\alpha_1$.

\begin{prop}\label{prop:feet_distance}
    For any $\pi/30>\epsilon>0$, there exists $R_0>0$ such that when $R>R_0$, the following holds: Suppose $\eta$ is a third connection for $\gamma_0$ such that the generated pants is $(R,\epsilon)$-cuff-good and good. Then $d^1_{S^{n-2}}(\eta'(A),-\eta'(B))<3\epsilon/2$.
\end{prop}

\begin{proof}
    The idea of proving this statement is similar to Part 4 of the proof of Theorem \ref{2.9}. Suppose $X_1,X_2$ and $Y$ are defined as above; then we know that $X_i^{-1}\phi(Y)$ is $\epsilon$-close to $\mathbf{e}$, and $X_2^{-1}\phi(X_1)$ is conjugate to $\Lambda$ which is also $\epsilon$-close to $\mathbf{e}$. Thus $X_1$ is $2\epsilon$-close to $X_2$, and $X_2$ is $\epsilon$-close to $\phi(X_1)$; hence $X_1$ is $3\epsilon$-close to $\phi(X_1)$. Then by Lemma \ref{2.7}, there exists $Q\in\mathcal{Q}$ such that $X_1Q$ is $3\epsilon/2$-close to $\mathbf{e}$ or $\mathbf{diag}(-1,-1,1,\ldots,1)$. Since the pants generated by $\gamma_0$ and $\eta$ is good, so we can rule out the second possibility, which means $X_1Q$ is $3\epsilon/2$-close to $\mathbf{e}$. We can assume that $X_1$ is $3\epsilon/2$-close to $\mathbf{e}$ by replacing $F$ with $F\cdot Q^{-1}$; here $Q$ is considered as an element in $\SO(n-2)$. Therefore $\phi(X_1)$ is also $3\epsilon/2$-close to $\mathbf{e}$, since $\phi(X_1)$ is conjugate to $X_1$ in $\mathrm{O}(n)$.
    
    On the other hand, we notice that the second vector of $E_1(x_1)$ is the parallel transport of $\eta'(A)$ to $x_1$, and the second vector of $F_1(x_1)\cdot R(\pi)$ is the parallel transport of $-\eta'(B)$ to $x_1$. Therefore by $\phi(X_1)$ sending $E_1(x_1)$ to $F_1(x_1)\cdot R(\pi)$ as $n$-frames and $d(\mathbf{e},\phi(X_1))<3\epsilon/2$, we have 
    \begin{equation*}
        d^1_{S^{n-2}}(\eta'(A),-\eta'(B))<3\epsilon/2.\qedhere
    \end{equation*}
\end{proof}

On the other hand, we want to show that when the two feet are sufficiently close, i.e., $d^1_{S^{n-2}}(\eta'(A),-\eta'(B))$ is sufficiently small, the distance between $X_1$ and $X_2$ is bounded by $3\epsilon/2$; this will ensure that we have enough third connections in the good region.

\begin{prop}\label{prop:monodromy_distance}
    For any $\pi/30>\epsilon>0$, there exists $R_0>0$ such that when $R>R_0$, the following holds: Suppose $\eta'(A), -\eta'(B)$ are unit normal vectors at $A$ and $B$ on $\gamma_0$ respectively (where we keep the notation even though we are not given $\eta$) such that $d^1_{S^{n-2}}(\eta'(A),-\eta'(B))<\epsilon/4$. Then
    \begin{equation}
        d(X_1,X_2)<\epsilon+2d^1_{S^{n-2}}(\eta'(A),-\eta'(B))<3\epsilon/2.
    \end{equation}
\end{prop}

\begin{rem}
    The change of the frames $E$ or $F$ will not affect the distance between $X_1$ and $X_2$; therefore, $d(X_1,X_2)$ is well-defined without choosing $E$ and $F$ and can be calculated from any choice of $E$ and $F$.
\end{rem}

\begin{proof}
    Let $\delta=d^1_{S^{n-2}}(\eta'(A),-\eta'(B))<\epsilon/4$. After bringing $\eta'(A)$ and $-\eta'(B)$ to $x_1$, let $P$ be a 2-plane in the normal space at $x_1$ to $\gamma_0$ that contains $\eta'(A)$ and $-\eta'(B)$. Then we can choose frames $E$ at $A$ and $F$ at $B$ such that $X_1$ is the rotation in $P$ that maps $\eta'(A)$ to $-\eta'(B)$ and has rotation angle $\delta$. So $d(\mathbf{e},X_1)=\delta$, and hence $d(\mathbf{e},\phi(X_1))=\delta$. Then by $d(\mathbf{e},X_2^{-1}\phi(X_1))<\epsilon$, we know that 
    \begin{equation}
        d(\mathbf{e},X_2)\leq d(\mathbf{e},\phi(X_1)) + d(\phi(X_1),X_2)  =  d(\mathbf{e},\phi(X_1)) + d(\mathbf{e},X_2^{-1}\phi(X_1))<\epsilon+\delta.
    \end{equation} 
    Thus we have
    \begin{equation*}
        d(X_1,X_2)\leq d(X_1,\mathbf{e})+ d(\mathbf{e},X_2)<\epsilon+2\delta<3\epsilon/2.\qedhere
    \end{equation*}
\end{proof}


What we really want to do is to coarsely estimate the number of good pants with boundary $\gamma_0$ that have their short feet in a given small region in $N^1(\gamma_0)$. In Section \ref{sec:ave-feet}, we introduce the \emph{average feet} of the long feet, and prove that these average feet are very close to the short feet. Then, in Section \ref{sec:volume}, for each prospective average foot, we consider the space of possible pairs of long feet with that average foot, and integrate the ``volume of the intersection'' $\Vol(B_\epsilon(X_1) \cap B_\epsilon(X_2))$ over this space of pairs of long feet. Our coarse estimate, which appears as Proposition \ref{prop:mu_Lebesgue}, will then follow from Propositions \ref{prop:feet_distance} and \ref{prop:monodromy_distance}. 
. 

\subsection{The average feet of good pants}\label{sec:ave-feet}

For a pants $P$ in a hyperbolic 2 or 3-manifold, a short foot on one cuff $\gamma$ is always a midpoint of the long feet on $\gamma$ in the unit normal bundle of $\gamma$, and $P$ is divided into two isometric right-angled hexagons. But the same properties do not hold in hyperbolic manifolds of dimension greater than 3. In this section, we will define the average of long feet for good pants in $M$ and estimate the difference between the new defined feet and the short feet in $N^1(\gamma)$.

Suppose that $\Pi$ is an oriented $(R,\epsilon)$-cuff-good and good pants with a boundary curve $\gamma_0$, for $\epsilon>0$ sufficiently small such that it satisfies all the above propositions, and $R$ sufficiently large which will be determined later. Let $\eta$ be the third connection for $\gamma_0$ with two long feet $u_1,u_2\in N^1(\gamma_0)$, and $n_1,n_2$ be the two short feet on $\gamma_0$ such that $u_i,n_i,u_{i+1}$ are ordered along the orientation of $\gamma_0$ for $i\in\Z/2\Z$.

Let us now define the \emph{average feet} $a_1, a_2$ on $\gamma_0$ of $\Pi$. There are two midpoints $b_1, b_2$ of the base points for $u_1$ and $u_2$, and we label them so that $u_i,b_i,u_{i+1}$ are ordered along the orientation of $\gamma_0$ (as with the $n_i$). Fixing $i$, we parallel transport $u_1$ and $u_2$ to $b_i$ (transporting each about one-quarter the length of $\gamma_0$).
By Proposition \ref{prop:feet_distance}, 
we know that the resulting unit normal vectors
are $3\epsilon/2$-close in the unit normal sphere at $b_i$. We then let $a_i$ be the midpoint of the short geodesic (in this unit sphere at $b_i$) connecting them.






Before stating the theorem, we state a simple result in hyperbolic geometry, and its proof is left to the reader.

\begin{lem}\label{distance d}
    There exists a universal constant $K_0,d_0>0$ such that when $d>d_0$ the following holds: Suppose $\gamma_1,\gamma_2$ are two geodesics in $\H^n$ such that the real distance between $\gamma_1$ and $\gamma_2$ is $d$. For any two points $X_1,X_2\in\gamma_2$, let $n_i\in N^1(\gamma_1)$ be the foot of the perpendicular geodesic from $X_i$ to $\gamma_1$, $i=1,2$. Then $n_1$ and $n_2$ are within distance $K_0e^{-d}$ in $N^1(\gamma_1)$.
\end{lem}


Now we show that the average feet are close to the short feet for good pants.

\begin{prop}\label{ave close to short}
    There exists a universal constant $K>0$ such that when $R$ is large enough, $a_i$ is $Ke^{-R}$-close to $n_i$, $i\in\Z/2\Z$.
\end{prop}

\begin{proof}
    Let $\gamma_1,\gamma_2$ be the other two cuffs of the pants $\Pi$ and $\eta_i$ be the short orthogeodesic connecting $\gamma_0$ and $\gamma_i$ for $i=1,2$, and $n_i$ is one foot of $\eta_i$. Let $\beta$ be the third connection for $\gamma_0$ with two long feet $u_1,u_2$. We then focus on the lifts of the cuffs and feet to $\H^{n}$ and $T\H^{n}$, and want to estimate the difference of lifts of $a_i$ and $n_i$.

    Let $\tilde{\gamma}_0$ be a lift of $\gamma_0$, and $\tilde{u}_1,\tilde{a}_1,\tilde{n}_1,\tilde{u}_2$ be the lifts of $u_1,a_1,n_1,u_2$ respectively such that $\tilde{a}_1$ and $\tilde{n}_1$ are the only lifts of $a_1$ and $n_1$, respectively, between $\tilde{u}_1$ and $\tilde{u}_2$.  Without loss of generality, we then only need to show that $\tilde{a}_1$ is $Ke^{-R}$-close to $\tilde{n}_1$ for some constant $K>0$.
    
    Let $\tilde{\beta}$ be the lift of $\beta$ such that $\tilde{u}_1$ is the foot for $\tilde{\beta}$, and $\tilde{\gamma}'_0$ be another lift of $\gamma_0$ such that $\tilde{\beta}$ connects $\tilde{\gamma}_0$ and $\tilde{\gamma}'_0$. Let $\tilde{\eta}_1$ and $\tilde{\gamma}_1$ be the lifts of $\eta_1$ and $\gamma_1$ such that $\tilde{n}_1$ is the foot for $\tilde{\eta}_1$ at $\tilde\gamma_0\cap\tilde\eta_1$ and $\tilde{\eta}_1$ is the orthogeodesic between $\tilde\gamma_0$ and $\tilde\gamma_1$. We then let $\tilde{\eta}'_1$ be the lift of $\eta_1$ connecting $\tilde{\gamma}_1$ and $\tilde{\gamma}'_0$, and $Q=\tilde{\eta}'_1\cap\tilde{\gamma}'_0$ and $Q_1=\tilde{\eta}'_1\cap\tilde{\gamma}_1$. We denote by $\partial\H^{n}$ the boundary at infinity of $\H^{n}$, and let $X_1,X_2$ be the endpoints of $\tilde{\gamma}_1$ in $\partial\H^{n}$ such that $Q_1$ is between $X_1$ and $\tilde{\gamma}_1\cap\tilde{\eta}_1$. We drop the perpendicular geodesics from $Q,Q_1,X_1$ and $X_2$ to $\tilde{\gamma}_0$ with foot $q,q_1,x_1$ and $x_2$ on $\tilde\gamma_0$, respectively.  Then we know that $\tilde{n}_1$ is the midpoint of $x_1$ and $x_2$ in $N^1(\tilde{\gamma}_0)$. On the other hand, we know that $\tilde{a}_1$ is the midpoint of $\tilde{u}_1$ and the foot of another lift of $\beta$ on $\tilde{\gamma}_0$. Hence it is sufficient to prove that $\tilde{u}_1$ is $O(e^{-R})$-close to $x_1$.
    \begin{figure}[htpb]
        \centering
        \includegraphics[width=0.8\textwidth]{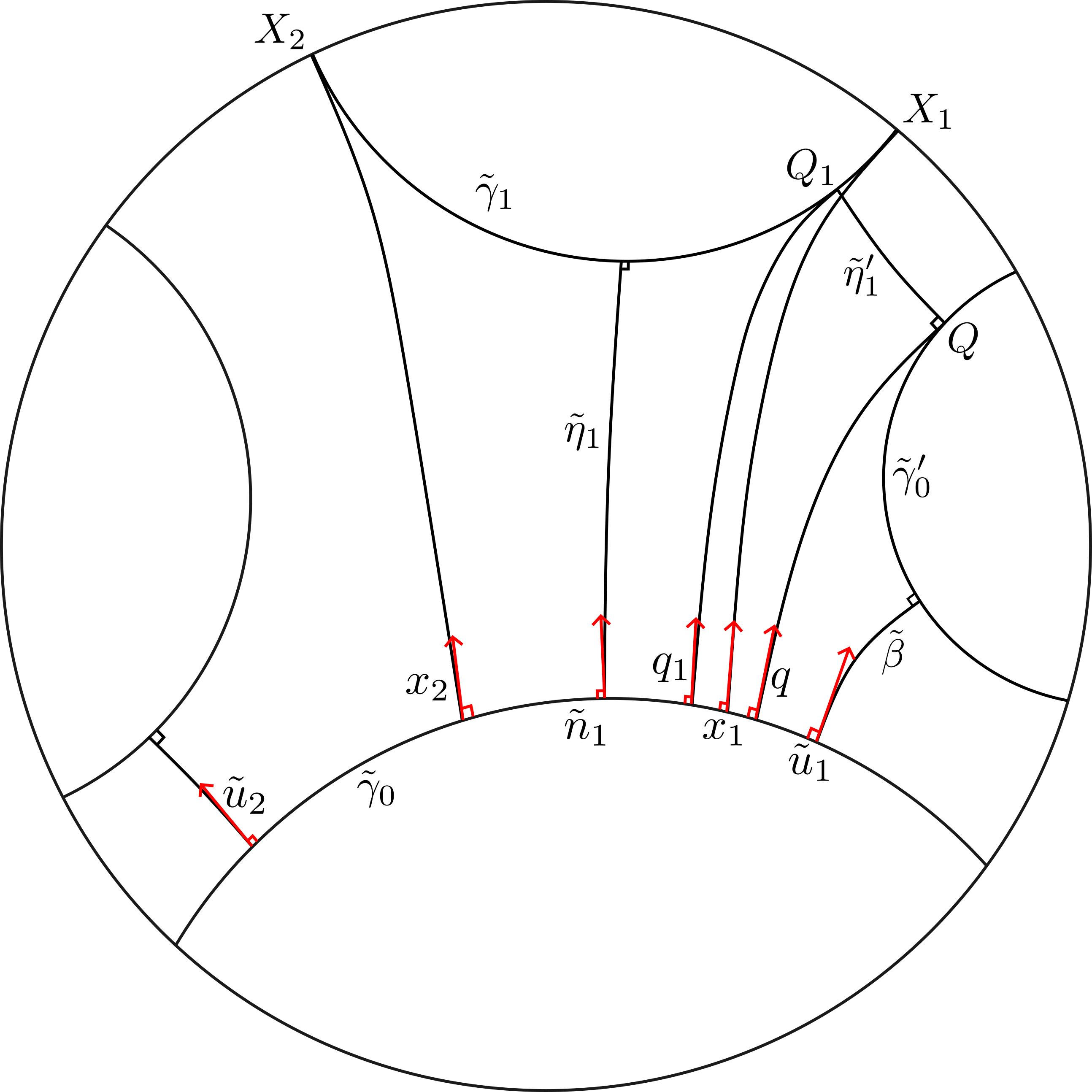}
        \caption{Lifts of cuffs and feet in $\H^{n}$}
        \label{fig:AverageFeet}
    \end{figure}

    We let $R>d_0$ be large enough, where $d_0$ is from Lemma \ref{distance d}. 
    Since $\gamma_1$ is a good curve, we have 
    \begin{equation*}
        2R-2\epsilon<d_{\H^{n}}(Q_1,\tilde{\eta}_1)=l(\gamma_1)<2R+2\epsilon.
    \end{equation*}
    We note that the distance between $\tilde{\gamma}_0$ and $\tilde{\gamma}_1$, which is the length of $\eta_1$, is bounded by a constant multiple of $e^{-R/2}$. Therefore by elementary hyperbolic geometry, $d_{\H^{n}}(Q_1,\tilde{\gamma}_0)>R$ and $q_1$ is $e^{-R}$-close to $x_1$. We notice that $d_{\H^{n}}(Q_1,Q)=l(\tilde{\eta}'_1)$ is also equal to the length of $\eta_1$, which is exponentially small. So by elementary hyperbolic geometry again and $d_{\H^{n}}(Q_1,\tilde{\gamma}_0)>R$, we know that $q$ is $e^{-R}$-close to $q_1$. Finally, by Lemma \ref{distance d}, $\tilde{u}_1$ is $K_0e^{-R}$-close to $q$. Hence $\tilde{u}_1$ is $(K_0+2)e^{-R}$-close to $x_1$.
\end{proof}

One may notice that the average feet can actually be defined without the given orientations of $\Pi$ and $\gamma$, and the above results all holds as well. Here the orientations are provided just for convenience. 

\subsection{The measure of feet of good pants} \label{sec:volume}
In this final subsection (of Section \ref{sec:counting}) we introduce the \emph{actual average foot measure} and \emph{estimated average foot measure}, and show that they are almost equivalent. We also show that the volume of intersection (as discussed in Section \ref{sec:computingmonodromy})
is bounded above and below, for all pairs of long feet that are close (after being brought together by parallel transport), which implies that the estimated average foot measure is quasi-uniform. This is our fundamental estimate, which we will use to prove the matching result in Section \ref{sec:matching_pants}.

In Section \ref{subsec:esti-foot}, 
we introduce the average foot measure on $N^1(\gamma_0)$,
along with the projection from the space of invariants to $N^1(\gamma_0)$.
We study the fibers of this projection and provide some estimates that we will use in Sections \ref{subsec:volume-intersection} and \ref{sec:measure-equiv}.
In Section \ref{subsec:volume-intersection},
we prove that the estimated average foot measure is bounded above and below; 
this is the new ingredient that allows us to go beyond \cite{Ham15}.
In Section \ref{sec:measure-equiv},
we relate the estimated average foot measure to the actual average foot measure, which is what we will use in Section \ref{sec:matching_pants} to match the good pants and build the surface.

\subsubsection{The average foot measures for $\gamma_0$} \label{subsec:esti-foot}
In this subsection we introduce a number of measures of two kinds: those that count actual third connections, and those that count predicted third connections, where the prediction comes, via exponential mixing, from \cite{KW19}. 
While we can think of third connections as living in a space of invariants, defined below, 
we want to push these measures down to a measure on $N^1(\gamma)$ by taking the average of the feet of the third connection. 

For each oriented closed geodesic $\gamma_0$, a framed oriented third connection $\eta$ can be described as a point $((u,(f_u,E)),(v,(f_v,F)),l,\Lambda)$ in
$$\mathbb{F}(\gamma_0)\equiv \mathcal{F}^\perp(\gamma_0)\times \mathcal{F}^\perp(\gamma_0)\times[0,\infty) \times \SO(n-1),$$
where $u,v$ are the two footpoints of $\eta$, $f_u$ is the tangent vector of $\eta$ at $u$, $f_v$ is the tangent vector of $\eta$ at $v$, and $E,F$ are orthonormal frames of the corresponding orthogonal complement in the normal spaces, $l$ is the length of $\eta$, and $\Lambda$ is the monodromy along $\eta$ (which lies in $\SO(n-1)$ because it preserves the tangent vector to $\eta$). We then have an action of $\SO(n-2)$ on each $\mathcal{F}^\perp(\gamma_0)$ factor of $\mathbb{F}(\gamma_0)$ as the rotation around the first vector of the frame (with the resulting left or right multiplication on $\SO(n-1)$), and therefore we have an action of $\SO(n-2)\times\SO(n-2)$ on $\mathbb{F}(\gamma_0)$ by 
\begin{equation}
    ((u,(f_u,E)),(v,(f_v,F)),l,\Lambda)\cdot(A,B)=((u,(f_u,E\cdot A)),(v,(f_v,F\cdot B)),l,A^{-1}\Lambda B).
\end{equation} 
Here $\SO(n-2)$ is regarded as a subgroup of $\SO(n-1)$ by $X\mapsto \begin{pmatrix}
    1 & \\
     & X
\end{pmatrix}.$
Let
\begin{equation}
    \I(\gamma_0)=\faktor{\mathbb{F}(\gamma_0)}{\SO(n-2)\times\SO(n-2)},
\end{equation}
which is a $\SO(n-1)$-fiber bundle over $N^1(\gamma_0)\times N^1(\gamma_0)\times [0,\infty)$; then $\I(\gamma_0)$ is the ``space of invariants'' for third connections for $\gamma_0$. We take the Riemannian metric on $\mathbb{F}(\gamma_0)$ to be the $L^2$-metric over each coordinate and push it down to $\I(\gamma_0)$. 
Locally, we let the measure on $\I(\gamma_0)$ be the product of the Lebesgue (or Haar) measures on $N^1(\gamma_0)$, $N^1(\gamma_0)$, and $\SO(n-1)$, and multiplied by $e^{2t}dt$ on $[0,\infty)$.


Let 
$\Pi_{R,\epsilon}(\gamma_0)$ be the set of pairs $(P,\eta)$ where $P$ is good and $(R,\epsilon)$-cuff-good in $M$ of which $\gamma_0$ is a boundary geodesic and $\eta$ is an oriented third connection of $P$ for $\gamma_0$, and define $\pi_0\colon\Pi(\gamma_0)\to\I(\gamma_0)$ by $\pi_0(P,\eta)=Q(\eta)$,
where $Q(\eta)$ is defined to be the image of $\eta$ in $\mathbb F (\gamma_0)$, pushed down to $\I(\gamma_0)$, by our above discussion.
Then $\Pi_{R,\epsilon}(\gamma_0)$ is a finite set and the counting measure on $\Pi_{R,\epsilon}(\gamma_0)$ generates a point measure $\tilde{\nu}$ on $\I(\gamma_0)$. We know that $\tilde{\nu}$ counts the third connections of all good and $(R,\epsilon)$-cuff good pants in $M$.

Now we want to define two \emph{good} regions to estimate $\tilde{\nu}$.
Let $\mathcal{R}(\gamma_0)\subset\I(\gamma_0)$ be the region
such that $(P, \eta)$ \emph{would} be good and $(R,\epsilon)$ cuff-good if $Q(\eta) \in \mathcal{R}(\gamma_0)\subset\I(\gamma_0)$. 
More formally,
we let $T_{\gamma_0}$ be the translation (defined up to conjugacy) in $\Isom(\H^{n})$ such that $\H^{n}/T_{\gamma_0}$ is the cover of $M$ associated to $\gamma_0$.
For any $\alpha\in\I(\gamma_0)$,
there is a unique pants group $P(\alpha)$ (which has $T_{\gamma_0}$ as one generator)
such that $\H^{n}/P(\alpha)$ has a third connection $\eta$ (from $\gamma_0$ to itself), 
such that $Q(\eta) = \alpha$.
We then say that $\alpha \in \mathcal{R}(\gamma_0)\subset\I(\gamma_0)$ if and only if $P(\alpha)$ is good and $(R, \epsilon)$-cuff-good. Essentially $\mathcal{R}(\gamma_0)$ corresponds to all potential good and $(R,\epsilon)$-cuff good pants with $\gamma_0$ as a cuff.

For any $\delta>0$, we want to define $\mathcal{R}_\delta(\gamma_0)\subset\I(\gamma_0)$ to be the region of third connections that satisfy our previous estimate in Section \ref{sec:computingmonodromy}. Formally, a point $\alpha\in\I(\gamma_0)$ is contained in $\mathcal{R}_\delta(\gamma_0)$ if and only if $P(\alpha)$ is good and there is a lift $((u,(f_u,E)),(v,(f_v,F)),l,\Lambda)$ of $\alpha$ in $\mathbb{F}(\gamma_0)$ such that:
\begin{enumerate}
    \item if $\eta$ is a third connection connecting $u$ and $v$ with length $l$, the lengths of two new cuffs are within $(2R-2\epsilon,2R+2\epsilon)$;
    \item suppose $X_1,X_2$ are the matrices of right actions defined as in Section \ref{sec:computingmonodromy} for frames at $u$ and $v$, where $f_u,f_v$ are regarded as $\eta'(A),-\eta'(B)$; then $\phi(Y)$ lies in the intersection of $\delta$-neighborhoods centered at $X_1$ and $X_2$.
    
\end{enumerate}

In the definitions of two good regions, the requirements about length are the same, but the assumptions for the monodromy are slightly different.

By definition, we have $\supp(\tilde{\nu})\subset \mathcal{R}(\gamma_0)$. By our observation (*) in Section \ref{sec:computingmonodromy}, there exists a universal constant $C_0>0$, such that
\begin{equation}\label{eq:first_good_region}
    \mathcal{R}_{\epsilon-C_0e^{-R}}(\gamma_0)
    \subset \mathcal{R}(\gamma_0)
    \subset\mathcal{R}_{\epsilon+C_0e^{-R}}(\gamma_0).
\end{equation}
By Theorem 4.9 in \cite{KW19}, we have the following proposition, which tells us that the measure $\tilde{\nu}$ of good third connections in $M$ can be approximated by $\mathbbm{1}_{\mathcal{R}(\gamma_0)}$, the measure of all potential good third connections.

\begin{prop}\label{original-counting}
    There exists $q=q(M)\leq1/2$ and a constant $C>0$ such that for any $\epsilon>0$, the following holds when $R$ is sufficiently large: Let $\gamma_0$ be an $(R,\epsilon)$-good curve and $\zeta=e^{-qR}$. Then for any $A\subset \mathcal{R}(\gamma_0)\subset \I(\gamma_0)$ we have
    \begin{equation}
        (1-\zeta)|\mathcal{N}_{-\zeta}(A)|\leq C\cdot\tilde{\nu}(A)\Vol(M) \leq(1+\zeta)|\mathcal{N}_\zeta(A)|,
    \end{equation}
    here $\mathcal{N}_{-\zeta}(A)=\I(\gamma_0)-\mathcal{N}_\zeta(\I(\gamma_0)-A)$.
\end{prop}

Since our plan is to match pants through their average feet, which are the ``midpoints" of long feet, as described in Section \ref{sec:ave-feet}, we want to know whether a unit normal vector could be the average foot of some pair of long feet. Hence we define the following functions to help us determine the distribution of average feet.

For $x,y \in N^1(\gamma_0)$ such that the distance between their basepoints is less than $l(\gamma_0)/2$, the point $m_x(y)\in N^1(\gamma_0)$ is derived by the following process:
\begin{enumerate}
    \item The basepoint of $m_x(y)$ is determined by extending the shortest geodesic segment along $\gamma_0$ from the basepoint of $y$ to the basepoint of $x$ by its length.
    \item Parallel transport $x,y$ along the extended geodesic segment in (1) to the basepoint of $m_x(y)$; we denote them by $x',y'$. Consider the geodesic from $y'$ to $x'$ in the unit normal sphere at the base point of $m_x(y)$; extend it by a distance equal to its length; the new endpoint is $m_x(y)$.
\end{enumerate}
Then $x$ can be viewed as \emph{a\ midpoint} of $y$ and $m_x(y)$ in the sense both of footpoints along the geodesic $\gamma_0$ and of unit vectors on the sphere after parallel transport. This function is essentially used to find another long foot given the average foot and one long foot.

Let $$\tilde{N}=\{(a,b):\text{the distance between the basepoints of }a\text{ and }b\text{ lies is }(l(\gamma_0)/4 - 1, l(\gamma_0)/4 + 1)\}$$
be a subset of $N^1(\gamma_0)\times N^1(\gamma_0)$. 
Then for $(a,b)\in \tilde{N}$, we let $\rho(a,b)=(m_a(b),b)\in N^1(\gamma_0)\times N^1(\gamma_0)$. 
Moreover, we can extend $\rho$ to a map 
from $\tilde{N}\times [0,\infty)$ to $N^1(\gamma_0)\times N^1(\gamma_0)\times[0,\infty)$ 
by preserving the third coordinate, let $\J(\gamma_0)=\rho^*(\I(\gamma_0))$ be the pullback bundle. 
For each point $(c,b,l,X)$ in our good region $\mathcal{R}(\gamma_0)$, there is only one preimage $(a,b,l,X)$ under $\rho_*$ such that:
\begin{enumerate}
    \item The oriented distance from $a$ to $b$ lies in $(l(\gamma_0)/4 - 1, l(\gamma_0)/4 + 1)$.
    \item After parallel transport to the same point on $\gamma_0$, the distance between $a$ and $b$ (on a unit normal sphere) is bounded by $\pi/4$.
\end{enumerate}
Here we need to assume that $\epsilon$ is no greater than some universal constant (by Proposition \ref{prop:feet_distance}) to satisfy the second requirement. Let ${\mathcal{S}}(\gamma_0)$ and ${\mathcal{S}}_\delta(\gamma_0)$ be the corresponding component of $(\rho_*)^{-1}(\mathcal{R}(\gamma_0))$ and $(\rho_*)^{-1}(\mathcal{R}_\delta(\gamma_0))$, respectively, that satisfies the above two conditions, and let $\hat{\nu}$ be the pullback measure of $\tilde{\nu}$ on $\J(\gamma_0)$. By \eqref{eq:first_good_region}, we have
\begin{equation}\label{eq:second_good_region}
    \mathcal{S}_{\epsilon-C_0e^{-R}}(\gamma_0)
    \subset \mathcal{S}(\gamma_0)
    \subset\mathcal{S}_{\epsilon+C_0e^{-R}}(\gamma_0).
\end{equation}
Moreover, since the Jacobian of $\rho$ is a constant, the result in Proposition \ref{original-counting} holds for $\hat{\nu}$ and ${\mathcal{S}}(\gamma_0)\subset\J(\gamma_0)$ as well, which indicates that $\hat{\nu}$ can be approximated by $\mathbbm{1}_{\mathcal{S}(\gamma_0)}$:
\begin{prop}\label{second-counting}
    There exists $q=q(M)\leq1/2$ and a constant $C>0$ such that for any $\epsilon>0$, the following holds when $R$ is sufficiently large: Let $\gamma_0$ be an $(R,\epsilon)$-good curve and $\zeta=e^{-qR}$. Then for any $A\subset {\mathcal{S}}(\gamma_0)\subset \J(\gamma_0)$ we have
    \begin{equation}
        (1-\zeta)|\mathcal{N}_{-\zeta}(A)|\leq C\cdot\hat{\nu}(A)\Vol(M) \leq(1+\zeta)|\mathcal{N}_\zeta(A)|,
    \end{equation}
    here $\mathcal{N}_{-\zeta}(A)=\J(\gamma_0)-\mathcal{N}_\zeta(\J(\gamma_0)-A)$.
\end{prop}

We now project the measures on the invariants of actual and potential pants down to the unit normal bundle $N^1(\gamma_0)$ to get the foot measures. Let $\pi\from \J(\gamma_0)\to N^1(\gamma_0)$ be the projection to the first $N^1(\gamma_0)$-coordinate and
$\nu_a=\pi_*(\hat{\nu})$. Then $\nu_a$ is the point measure on $N^1(\gamma_0)$ of the average feet of all good and $(R,\epsilon)$-cuff-good pants in $M$, and we call $\nu_a$ the \emph{actual average foot measure}. Let $\mu_a=\pi_*(\mathbbm{1}_{\mathcal{S}(\gamma_0)})$; it is the measure of average feet of all potential good and $(R,\epsilon)$-good pants. We call $\mu_a$ the \emph{estimated average foot measure}.
The measures $\hat{\nu},\nu_a$ and $\mu_a$ actually depend on the good curve $\gamma_0$, and we will keep using the notations when there is no ambiguity. 
Clearly,
the estimated average foot measure $\mu_a$ is invariant by the centralizer for the element of $\SO(n, 1)$ associated to $\gamma_0$, including translation along $\gamma_0$ and performing the antipodal map on each unit normal spheres.
We expect to push Proposition \ref{second-counting} down to the unit normal bundle and have a similar estimate for $\nu_a$ and $\mu_a$, meaning that $\nu_a$ can be approximated by $\mu_a$.
The result will be stated in Section \ref{sec:measure-equiv}, and we will prove it through fiberwise volume estimates for $\tilde{\nu}$ and $\mathbbm{1}_{\mathcal{S}(\gamma_0)}$.



We will denote by $v$ an arbitrary unit normal vector in $N^1(\gamma_0)$ from now on for our fiberwise estimate, and let $F^v=\pi^{-1}(v)$ be the preimage of $v$ in $\J(\gamma_0)$; then $F^v$ is a $\SO(n-1)$-bundle over $N^1(\gamma_0)\times[0,\infty)$.
    
For the intersection of each fiber $F^v$ with the good regions, we can actually think of it as a subset of a product space, where we can estimate the volumes of the spaces of lengths and monodromies separately. By Lemma \ref{lem:third-conn-length}, we know that if $(v,y,l,X)\in F^v\cap \mathcal{S}_\delta(\gamma_0)$, then $l\in (R+2\ln2-4\epsilon, R+2\ln2+4\epsilon)$ and the distance from $v$ to $y$ is $4\epsilon$-close to $R/2$ (here the distance from $v$ to $y$ is half of the distance between two footpoints of a \emph{potential} third connection).
By the parallel transport along $\gamma_0$ (of distance about $R/2$) from the basepoint of $y$ to that of $v$, we can identify the unit normal spheres at these two points. Then $y$ can be described as a point in the unit normal sphere at the basepoint of $v$ together with the distance from the basepoint of $v$ to that of $y$.
Since the restriction of $N^1(\gamma_0)$ to an arc of $\gamma_0$ is a trivial bundle, then for $(v,y,l,X)\in F^v\cap \mathcal{S}_\delta(\gamma_0)$, its second coordinate can viewed as in the product space $S^{n-2}\times (R/2-1,R/2+1)$. Hence we can think of $F^v\cap \mathcal{S}_{\delta}(\gamma_0)$ as a subset of a $\SO(n-1)$-bundle over $$S^{n-2}\times (R/2-1,R/2+1)\times [0,\infty).$$
Moreover, by our previous discussion in Section \ref{sec:computingmonodromy}, the fiber of this bundle only depends on its projection on $S^{n-2}$. Thus this bundle is a product space of $D^v$ and $(R/2-1,R/2+1)\times [0,\infty)$, where $D^v$ is a $\SO(n-1)$-bundle over $S^{n-2}$.
Let $\hat{\mathcal{S}}^v_\delta$ be the projection of $F^v\cap \mathcal{S}_{\delta}(\gamma_0)$ on $D^v$. Then by Lemma \ref{lem:third-conn-length}, 
\begin{equation}\label{eq:fibernbhd1}
    \hat{\mathcal{S}}^v_\delta\times \hat{L}_{R,\epsilon-Ce^{-R}}
    \subset F^v\cap \mathcal{S}_{\delta}(\gamma_0) 
    \subset \hat{\mathcal{S}}^v_\delta\times \hat{L}_{R,\epsilon+Ce^{-R}}.
\end{equation}
Here $\hat{L}_{R,\delta}\subset\R\times\R$ with measure $e^{2y}\,dx\,dy$ is the set of all pairs $(x,y)$ satisfying:
\begin{enumerate}
    \item $|2x+y-2\ln2-2R|<2\delta$;
    \item $|(l(\gamma_0)-2x)+y-2\ln2-2R|<2\delta$.
\end{enumerate}
Similar to \eqref{eq:diamondregion1}, we have 
\begin{equation}\label{eq:hat_L}
    \Vol(\hat{L}_{R,\delta})=4e^{4R-l(\gamma_0)}(e^{2\delta}-e^{-2\delta})^2
    \sim 64 \delta^2 e^{4R-l(\gamma_0)},
\end{equation}
when $\delta\to0$. 

\subsubsection{The volume of the space of good monodromies}\label{subsec:volume-intersection}
In this subsection, we give an estimate of the volume of $\hat{\mathcal{S}}^v_\delta$ for any $\delta$ that is close to $\epsilon$ and any fixed $v\in N^1(\gamma_0)$.
We first present a simple and rough estimate of the volume of the intersection of two $\epsilon$-balls in $\SO(n-1)$ in this subsection.

We first introduce some notations that can simplify our computation:
\begin{itemize}
    \item By $A\asymp_n B$, we mean that there exists a constant $K=K(n)$ such that 
    \begin{equation}
        \frac{B}{K(n)}<A<K(n)B.
    \end{equation}

    \item By $A\gtfrown_n B$, we mean that there exists a constant $K=K(n)$ such that 
    \begin{equation}
        A>K(n)B.
    \end{equation}

    \item By $A\ltfrown_n B$, we mean that there exists a constant $K=K(n)$ such that 
    \begin{equation}
        A<K(n)B.
    \end{equation}
\end{itemize}
Then $A\asymp_n B$ if and only if $A\gtfrown_n B$ and $A\ltfrown_n B$.

For $X \in \SO(n-1)$, let $V(X,r)$ denote the volume of the intersection of two balls of radius $r$ in $\SO(n-1)$, whose centers are at $\mathbf{e}$  and $X$. 
We let $V_m$ denote the volume of a ball of unit radius in $\R^m$, and we let $k=(n-1)(n-2)/2$ be the dimension of $\SO(n-1)$. When there is definite overlap to the balls we have a simple estimate:
\begin{prop}
    If $\norm{X} < 5r/3$, and $r$ is sufficiently small (given $n$),
    then 
\begin{equation}\label{lower bound}
    V(X,r)\asymp_n r^{k}.
\end{equation}
\end{prop}
\begin{proof}
For $\norm{X}$ sufficiently small,
we can find $\sqrt{X}\in \SO(n-1)$ such that $\sqrt{X}^2 = X$ and $\norm{\sqrt{X}} = \frac12 \norm{X}$, and $\sqrt{X}$ is the midpoint of the minimal geodesic from $\mathbf{e}$ to $X$. 
Then $B_{\frac r6}(\sqrt X) \subset B_r(\mathbf{e}) \cap B_r(X)$,
and 
$$
\Vol(B_{\frac r6}(\sqrt X)) = \Vol(B_{\frac{r}{6}}(\mathbf{e})) \ge C_n r^k,
$$
for $r$ sufficiently small and a universal constant $C_n$ depending on $n$,
by the Bishop-Gromov inequality.
\end{proof}

When $\norm{X}$ is close to $2r$, we have the following upper bound on $V(X, r)$.
\begin{prop} \label{prop:nearly disjoint}
    Let $\kappa = 2r - \norm{X}$, and suppose that $r < \epsilon_n$ and $\kappa < r$.
    Then 
    \begin{equation} \label{eq:nearly disjoint}
        V(X, r) < 2^k V_{k-1} \kappa^{\frac{k+1}2}r^{\frac{k-1}2},
    \end{equation}
    so $V(X, r) \ltfrown_n \kappa^{\frac{k+1}2}r^{\frac{k-1}2} \le \kappa^2 r^{k-2}$ when $n\geq 4$.
\end{prop}

This proposition, which, like the previous one, is more about distance and volume with control on curvature, follows from the following two theorems, whose proofs are relatively standard and outlined only briefly. If $x$ and $y$ are points in a complete Riemannian manifold, with a unique minimal geodesic between them, we let $\geod xy$ denote that geodesic segment. 
\begin{thm} \label{thm:triangle}
Suppose that $M^m$ is a complete connected Riemannian manifold, with all sectional curvatures bounded above and below by a certain universal constant $\epsilon_0 > 0$, and injectivity radius at least 4. Suppose that $x, y, z \in M$, with $d(x, z), d(y, z)  \le 1$ and $d(x, y) = 2 - s$. Let $w$ be the midpoint of $x$ and $y$ along $\geod x y$. 
Then there is a point $\hat z \in \geod x y$ such that
\begin{enumerate}
    \item 
    $d(w, \hat z) < s$ (where the distance is computed along $\alpha$),
    \item 
    $\geod z {\hat z}$ meets $\geod x y$ at right angles, and 
    \item 
    $d(z, \hat z) < 2\sqrt{s}$.
\end{enumerate}
\end{thm}
\begin{proof}
    The Theorem follows in $\R^m$, even when we divide the upper bounds in the inequalities by 2. 
    Applying Toponogov's Theorem, 
    we can then prove it in $M^m$ with sufficiently small bounds on sectional curvatures.
\end{proof}

For any geodesic segment $\alpha$ in a complete Riemannian manifold $M^m$, we let $B_r(\alpha)$ be all points in $M$ that can be reached by geodesic flow from a unit normal vector for $\alpha$ by time at most $r$; then $B_r(\alpha)$ is a cylinder-like region around $\alpha$. We let $|\alpha|$ denote the length of $\alpha$. 
\begin{thm} \label{thm:cylinder}
    Suppose that $M^m$ has non-negative sectional curvatures. Then $\Vol(B_r(\alpha)) \le V_{m-1} r^{m-1} |\alpha|$. 
\end{thm}
\begin{proof}
    This follows from Gauss's lemma and the Rauch comparison theorem.
\end{proof}

\begin{proof}[Proof of Proposition \ref{prop:nearly disjoint}]
Suppose $r$ is sufficiently small.
By Theorem \ref{thm:triangle}, 
applied after rescaling the metric on $\SO(n-1)$ by $1/r$,
for every $z \in B_r(\mathbf{e}) \cap B_r(X)$, 
we can find $\hat z$ on $\geod {\mathbf{e}}X$ within $\kappa$ of $\sqrt X$, 
such that $\geod {\hat z} z$ is orthogonal to $\geod {\mathbf{e}}X$
and $|\geod {\hat z} z| < 2 \sqrt{\kappa}$.
The Proposition then follows from Theorem \ref{thm:cylinder}, since $\SO(n-1)$ has non-negative sectional curvatures.
\end{proof}

Suppose $v\in N^1(\gamma_0)$;  for $a,b\in N^1(\gamma_0)$ that have the same basepoint as $v$, we arbitrarily choose two collections $E,F$ of unit vectors such that $(a,E)$ and $(b,F)$ are orthonormal frames of the normal space to $\gamma_0$ at the basepoint of $t$ with opposite orientations. Let $X_1,X_2\in\SO(n-1)$ be such that
\begin{equation}
    (-a,E)\cdot X_1=(b,F)
\end{equation}
and
\begin{equation}
    \Lambda\cdot(a,E)\cdot X_2=(-b,F),
\end{equation}
where $\Lambda$ is the monodromy of $\gamma_0$; we then have that $\Lambda$ is conjugate to $X_2^{-1}\phi(X_1)$ (similar to the computation after Lemma \ref{2.3}). We define $W_v(a,b)=X_1^{-1} X_2$ and $f_v(a,b)=d(X_1,X_2)$; then $f_v$ is well-defined and $W_v$ is well-defined up to conjugacy, independent of the choices of frames $E$ and $F$.



We recall that $m_{x}(y)\in N^1(\gamma_0)$ is defined for $x$ and $y$ when the distance between the basepoints of $x$ and $y$ is less than $l(\gamma_0)/2$. For simplicity, we will slightly abuse notation by using $m$ to denote the same kind of maps for two points on $S^{n-2}$. 

By our definition of $\hat{\S}_\delta^v\subset D^v$, the projection of $\hat{\S}_\delta^v$ on $S^{n-2}$ (as the base space of $D^v$) consists of $y\in S^{n-2}$ such that $f_v(y,m_v(y))<2\delta$. Here $v$ and $y$ are considered as points in the same $S^{n-2}$ by parallel transport, since $v$ is fixed.
Therefore for we have
\begin{align*}
    \Vol(\hat{\S}_\delta^v)
    &=\int_{S^{n-2}} V(W_v(y,m_v(y)),\delta)\, d\sigma.
\end{align*}
We now can present an estimate of this volume for $\delta$ that is close to $\epsilon$.

\begin{prop}\label{mu estimate}
  For any $v\in N^1(\gamma_0)$ and $\delta \in (8\epsilon/9, 7\epsilon/6)$, we have
 \begin{equation}
       \Vol(\hat{\S}_\delta^v) \asymp_n\epsilon^{(n-2)(n+1)/2}.
    \end{equation}\end{prop}

\begin{proof}
    It is equivalent to show $$\Vol(\hat{\S}_\delta^v)\ltfrown_n\epsilon^{(n-2)(n+1)/2}$$ and $$\Vol(\hat{\S}_\delta^v)\gtfrown_n\epsilon^{(n-2)(n+1)/2}.$$ 
    
    By the similar reasoning in the proof to Proposition \ref{prop:feet_distance}, we know that for any $v\in N^1(\gamma_0)$, if $f_v(y,m_v(y))<2\delta$, then
    \begin{equation}
        d_{S^{n-2}}(y,m_v(y))<\frac{\epsilon+2\delta}{2}<\frac{5\epsilon}{3}.
    \end{equation}
    We then have
    \begin{equation}\label{eq:upperbound}
        \begin{aligned}
            \Vol(\hat{\S}_\delta^v)&=\int_{S^{n-2}} V(W_v(y,m_v(y)),\delta) \,d\sigma \\
            &\leq \int_{S^{n-2}} \mathbbm{1}_{\{y:d_{S^{n-2}}(y,m_v(y))<5\epsilon/3\}} V(W_v(y,m_v(y)),\delta) \, d\sigma \\
            &\asymp_n \int_{S^{n-2}} \mathbbm{1}_{\{y:d_{S^{n-2}}(y,v)<3\delta/4\}} \delta^{k} \, d\sigma \\
            &\asymp_n \Vol({\{y:d_{S^{n-2}}(y,v)<3\delta/4\}})\cdot \epsilon^{k}\\
            &\asymp_n \epsilon^{n-2+k}=\epsilon^{(n-2)(n+1)/2}.
        \end{aligned}
    \end{equation}
    Hence $\Vol(\hat{\S}_\delta^v)\ltfrown_n \epsilon^{(n-2)(n+1)/2}$.

    By Proposition \ref{prop:monodromy_distance}, for $x,y$ on the unit normal sphere, if $d_{S^{n-2}}(x,y)<3\delta/16$, then 
    \begin{equation}
        f_v(x,y)<\epsilon+2d_{S^{n-2}}(x,y)\leq 9\delta/8 + 3\delta/8=3\delta/2<5\delta/3.
    \end{equation}
    Together with \eqref{lower bound}, we have
    \begin{equation}\label{eq:lowerbound}
        \begin{aligned}
            \Vol(\hat{\S}_\delta^v)&=\int_{S^{n-2}} V(W_v(y,m_v(y)),\delta) \,d\sigma \\
            &\geq \int_{S^{n-2}} \mathbbm{1}_{\{y:f_v(y,m_v(y))<3\delta/2\}}  V(W_v(y,m_v(y)),\delta) \,d\sigma \\
            &\asymp_n \int_{S^{n-2}} \mathbbm{1}_{\{y:f_v(y,m_v(y))<3\delta/2\}} \cdot \delta^{k} \,d\sigma \\
            &\geq \int_{S^{n-2}} \mathbbm{1}_{\{y:d_{S^{n-2}}(y,m_v(y))<3\delta/16\}}\cdot \delta^{k}\, d\sigma\\
            &= \int_{S^{n-2}} \mathbbm{1}_{\{y:d_{S^{n-2}}(y,v)<3\delta/32\}}\cdot \delta^{k}\, d\sigma\\
            &=\Vol({\{y:d_{S^{n-2}}(y,v)<3\delta/32\}})\cdot \delta^{k}\\
            &\asymp_n \delta^{(n-2)(n+1)/2)}\asymp_n \epsilon^{(n-2)(n+1)/2}.
        \end{aligned}
    \end{equation}
    Therefore $\Vol(\hat{\S}_\delta^v)\gtfrown_n\epsilon^{(n-2)(n+1)/2}$.\qedhere
 \end{proof}

With the above proposition we can prove the main result of  this section, that the  the estimated average foot measure $\mu_a$ is bounded above and below by the Lebesgue measure; we also call that $\mu_a$ is \emph{quasi-uniform} with respect to the Lebesgue measure. 

\begin{prop}\label{prop:mu_Lebesgue}
    There exists $B_0,\epsilon_0>0$ such that for any $\epsilon\in(0,\epsilon_0)$, the following holds when $R$ is sufficiently large: Suppose $\gamma_0$ is an $(R,\epsilon)$-good curve, and $\lambda$ is the Lebesgue measure on $N^1(\gamma_0)$. Then for any $v\in N^1(\gamma_0)$,
    \begin{equation}
       B_0^{-1}e^{2R} \epsilon^{(n^2-n+2)/2} 
       < \frac{d\mu_a}{d\lambda}\Big\vert_v 
       < B_0 e^{2R} \epsilon^{(n^2-n+2)/2}.
    \end{equation}
\end{prop}

\begin{proof}
    By definition, we know that for any $v\in N^1(\gamma_0)$,
    \begin{equation}
        \frac{d\mu_a}{d\lambda}\Big\vert_v=|F^v\cap \mathcal{S}(\gamma_0)|.
    \end{equation}
    Then it suffices to show $|F^v\cap \S(\gamma_0)|\asymp_n e^{2R} \epsilon^{(n^2-n+2)/2}$.
    
    By \eqref{eq:second_good_region}, \eqref{eq:fibernbhd1}, \eqref{eq:hat_L}, and Proposition \ref{mu estimate}, we have
    \begin{equation}
    \begin{aligned}
            |F^v\cap \S(\gamma_0)|&\leq |F^v\cap \S_{\epsilon+C_0e^{-R}}(\gamma_0)|\\
            &\leq |\hat{\S}^v_{\epsilon+C_0e^{-R}} \times \hat{L}_{R,\epsilon+Ce^{-R}}|\\
            &=\Vol(\hat{\S}^v_{\epsilon+C_0e^{-R}}) \times \Vol(\hat{L}_{R,\epsilon+Ce^{-R}})\\
            &\ltfrown_n \epsilon^{(n-2)(n+1)/2} 
             \times 800 e^{2R} \epsilon^2
            \asymp_n e^{2R} \epsilon^{(n^2-n+2)/2}.
        \end{aligned}
    \end{equation}
Similarly, we have $|F^v\cap \S(\gamma_0)|\gtfrown_n e^{2R} \epsilon^{(n^2-n+2)/2}$; then the result is proved.
\end{proof}

\subsubsection{The actual and the estimated average foot measures}\label{sec:measure-equiv}

In this subsection, we compare the actual average foot measure to the estimated, to have a similar result to Proposition \ref{second-counting}. 
We start with a lemma for the fibers of good regions in $\J(\gamma_0)$ over a point in $N^1(\gamma_0)$.

\begin{lem}\label{lem:fiberwisenbhd}
    There exists $C_1>0$ such that when $R$ is large enough: for any $v\in N^1(\gamma_0)$,
    \begin{equation}
        F^v\cap N_{\zeta}(\mathcal{S}(\gamma_0))
        \subset F^v\cap \mathcal{S}_{\epsilon+C_1\zeta}(\gamma_0),
    \end{equation}
    where $\zeta=e^{-qR}$ as in Proposition \ref{second-counting}.
\end{lem}

\begin{proof}
    Fix $v\in N^1(\gamma_0)$. 
    By \eqref{eq:second_good_region}, we have
    \begin{equation}\label{eq:F^t_nbhd1}
        F^v\cap N_{\zeta}(\mathcal{S}(\gamma_0))\subset 
        F^v\cap N_{\zeta}(\mathcal{S}_{\epsilon+C_0e^{-R}}(\gamma_0)).
    \end{equation}
    Suppose $\alpha=(v,y,l,Y)\in F^v\cap N_{\zeta}(\mathcal{S}_{\epsilon+C_0e^{-R}}(\gamma_0))$ for some $y\in N^1(\gamma_0)$, $l>0$ and $Y$ in the $\SO(n-1)$-fiber at $(v,y)$. Then there exists $\beta=(v',y',l',Y')\in\mathcal{S}_{\epsilon+C_0e^{-R}}(\gamma_0)$ such that $d(\alpha,\beta)\leq\zeta$. 
    Let $x=m_v(y)$ and $x'=m_v'(y')$ in $N^1(\gamma_0)$; here the distance between the basepoints of $v$ and $y$ is about $R/2$, so $m_v(y)$ is well-defined, and so is $m_v'(y')$. Since $\rho$ is a Lipschitz map, there exists a universal constant $C_2$, such that $(x,y,l,Y)$ is $C_2\zeta$-close to $(x',y',l',Y')$ in $\I(\gamma_0)$. Therefore we can find the lifts $(E_1,E_2,l,Y)$ and $(E'_1,E'_2,l',Y')$ of $(x,y,l,Y)$ and $(x',y',l',Y')$ respectively, such that 
    \begin{enumerate}
        \item $(E_1,E_2,l,Y)$ is $C_2\zeta$-close to $(E'_1,E'_2,l',Y')$ in $\mathbb{F}(\gamma_0)$;
        \item $\phi(Y')$ is in the $(\epsilon+C_0e^{-R})$-neighborhood of $X'_1$ and $X'_2$.
    \end{enumerate}
    Here $X'_i$ (and $X_i$) are the matrices for $(E'_1,E'_2,l',Y')$ (and $(E_1,E_2,l,Y)$) as described in Condition 2 in the definition of $\mathcal{R}_\delta(\gamma_0)$. Now since $d(E_i,E'_i)<C_2\zeta$, so 
    \begin{equation}
        d(X_1,X'_1),d(X_2,X'_2)<2C_2\zeta.
    \end{equation}
    Then by $d(Y,Y')<C_2\zeta$, we know $Y$ is in the $(\epsilon+C_0e^{-R}+3C_2\zeta)$-neighborhood of $X_1$ and $X_2$. Let $C_1=C_0+3C_2$; then by $q\leq1$, we have $\alpha\in F^v\cap \mathcal{S}_{\epsilon+C_0e^{-R}+3C_2\zeta}(\gamma_0)\subset F^v\cap \mathcal{S}_{\epsilon+C_1\zeta}(\gamma_0)$, and by \eqref{eq:F^t_nbhd1},
    \begin{equation*}
        F^v\cap N_{\zeta}(\mathcal{S}(\gamma_0))\subset 
        F^v\cap N_{\zeta}(\mathcal{S}_{\epsilon+C_0e^{-R}}(\gamma_0))
        \subset F^v\cap \mathcal{S}_{\epsilon+C_1\zeta}(\gamma_0).\qedhere
    \end{equation*}
        
\end{proof}

We now compare the actual average foot measure and the estimated average foot measure on $N^1(\gamma_0)$.

\begin{prop}\label{maincounting}
    There exists $q=q(M)\leq1/2$ and constant $C,L>0$ such that for any $\epsilon>0$, the following holds when $R$ is sufficiently large: Let $\gamma_0$ be an $(R,\epsilon)$-good curve and $\zeta=e^{-qR}$. Then for any $B\subset N^1(\gamma_0)$, we have
    \begin{equation}
        (1-L\zeta/\epsilon)\mu_a(\mathcal{N}_{-\zeta}(B))\leq C\cdot\nu_a(B)\Vol(M) \leq(1+L\zeta/\epsilon)\mu_a(\mathcal{N}_{\zeta}(B)),
    \end{equation}
    here $\mathcal{N}_{-\zeta}(B)=N^1(\gamma_0)-\mathcal{N}_{\zeta}(N^1(\gamma_0)-B)$.
\end{prop}

\begin{proof}
    For any $B\subset N^1(\gamma_0)$, we have $\nu_a(B)=\hat{\nu}(\pi^{-1}(B))$. Then by Proposition \ref{second-counting}, we have
    \begin{equation}
        \begin{aligned}
            &(1-\zeta)|\mathcal{N}_{-\zeta}(\pi^{-1}(B)\cap \mathcal{S}(\gamma_0))|\\
            \leq & C_0\cdot\hat{\nu}(\pi^{-1}(B))\cdot\Vol(M)
            =  C_0 \cdot \nu_a(B) \cdot \Vol(M)\\
            \leq & (1+\zeta)|\mathcal{N}_{\zeta}(\pi^{-1}(B)\cap \mathcal{S}(\gamma_0))|,
        \end{aligned}
    \end{equation}
    here $C_0$ is the constant $C$ in Proposition \ref{second-counting}. Thus we only need to show that there exists $L>0$ such that 
    \begin{equation}
        |\mathcal{N}_{\zeta}(\pi^{-1}(B)\cap \mathcal{S}(\gamma_0))|<(1+L\zeta/\epsilon)\mu_a(\mathcal{N}_{\zeta}(B)),
    \end{equation}
    and the proof of the other side of the inequality will be the same.

    We notice that 
    \begin{equation}
        \begin{aligned}
            \mathcal{N}_{\zeta}(\pi^{-1}(B)\cap \mathcal{S}(\gamma_0)) \subset \mathcal{N}_{\zeta}(\pi^{-1}(B)) \cap \mathcal{N}_{\zeta}(\mathcal{S}(\gamma_0))=\pi^{-1}(\mathcal{N}_{\zeta}(B)) \cap \mathcal{N}_{\zeta}(\mathcal{S}(\gamma_0)),
        \end{aligned}
    \end{equation}
    and
    \begin{equation}
        \mu_a(\mathcal{N}_{\zeta}(B))=|\pi^{-1}(\mathcal{N}_{\zeta}(B))\cap \mathcal{S}(\gamma_0)|.
    \end{equation}
    Therefore it suffices to show that there exists $L>0$ such that
    \begin{equation}
        |\pi^{-1}(\mathcal{N}_{\zeta}(B)) \cap \mathcal{N}_{\zeta}(\mathcal{S}(\gamma_0))|\leq (1+L\zeta/\epsilon) |\pi^{-1}(\mathcal{N}_{\zeta}(B))\cap \mathcal{S}(\gamma_0)|.
    \end{equation}
    We will prove this inequality fiberwise, i.e., we want to prove that there exists $L>0$ such that for any $v\in N^1(\gamma_0)$,
    \begin{equation}\label{eq:fiberwise-inequality}
        |F^v \cap \mathcal{N}_{\zeta}(\mathcal{S}(\gamma_0))|\leq (1+L\zeta/\epsilon) |F^v \cap \mathcal{S}(\gamma_0)|.
    \end{equation}

    By Lemma \ref{lem:fiberwisenbhd} and \eqref{eq:fibernbhd1}, we have
    \begin{equation}\label{eq:leftside}
        F^v \cap \mathcal{N}_{\zeta}(\mathcal{S}(\gamma_0))\subset F^v \cap \mathcal{S}_{\epsilon+C_1\zeta}(\gamma_0) \subset \hat{\mathcal{S}}_{\epsilon+C_1\zeta}^v \times \hat{L}_{R,\epsilon+Ce^{-R}}.
    \end{equation}
    By \eqref{eq:second_good_region} and \eqref{eq:fibernbhd1}, we have
    \begin{equation}\label{eq:rightside}
        F^v \cap \mathcal{S}(\gamma_0) \supset F^v \cap \mathcal{S}_{\epsilon-C_0e^{-R}}\supset \hat{\mathcal{S}}_{\epsilon-C_0e^{-R}}^v \times \hat{L}_{R,\epsilon-Ce^{-R}}.
    \end{equation}
    By \eqref{eq:hat_L} and $\zeta=e^{-qR}$ for $q\leq1/2$, we know that when $R$ is sufficiently large,
    \begin{equation}\label{eq:L-estimate}
        \Vol(\hat{L}_{R,\epsilon+Ce^{-R}})\leq (1+\zeta)\Vol(\hat{L}_{R,\epsilon-Ce^{-R}}).
    \end{equation}
    Hence now we only need to compare $\Vol(\hat{\mathcal{S}}_{\epsilon+C_1\zeta}^v)$ and $\Vol(\hat{\mathcal{S}}_{\epsilon-C_0e^{-R}}^v)$.

    Let $b\from D^v\to S^{n-2}$ be the projection map, and for any $s\in S^{n-2}$, let $G^s=b^{-1}(s)\cong \SO(n-1)$. By definition, we have $\hat{\mathcal{S}}_{\epsilon-C_0e^{-R}}^v\subset\hat{\mathcal{S}}_{\epsilon+C_1\zeta}^v$. Suppose $u\in b(\hat{\mathcal{S}}_{\epsilon+C_1\zeta}^v)\subset S^{n-2}$; then $f_v(u,m_v(u))\leq 2(\epsilon+C_1\zeta)$.

    We now want to obtain an upper bound on $\Vol(G^u\cap(\hat{\mathcal{S}}_{\epsilon+C_1\zeta}^v-\hat{\mathcal{S}}_{\epsilon-C_0e^{-R}}^v))$. We have two cases, depending on $f_v(u, m_v(u))$:
    \begin{itemize}
        \item If $2(\epsilon-C_0e^{-R})<f_v(u,m_v(u))\leq 2(\epsilon+C_1\zeta)$, then $u\notin b(\hat{\mathcal{S}}_{\epsilon-C_0e^{-R}}^v)$. Thus Proposition \ref{prop:nearly disjoint}, when $R$ is sufficiently large, we have
        \begin{align*}\label{eq:case1}
                \Vol(G^u\cap(\hat{\mathcal{S}}_{\epsilon+C_1\zeta}^v-\hat{\mathcal{S}}_{\epsilon-C_0e^{-R}}^v))&=\Vol(G^u\cap\hat{\mathcal{S}}_{\epsilon+C_1\zeta}^v)\\
                & =V(W_v(u,m_v(u)),\epsilon+C_1\zeta)\\
                &\ltfrown_n (\epsilon+C_1\zeta)^{\frac{k-1}{2}} (2\epsilon+2C_1\zeta-f_v(u,m_v(u)))^{\frac{k+1}{2}}\\
                &< (2\epsilon)^{\frac{k-1}{2}}(2C_1\zeta + 2C_0 e^{-R})^{\frac{k+1}{2}}\\
                &\ltfrown_n \epsilon^{\frac{k-1}{2}}\zeta^{\frac{k+1}{2}}. \numberthis
        \end{align*}

        \item If $f_v(u,m_v(u))\leq 2(\epsilon-C_0e^{-R})$, then $u\in b(\hat{\mathcal{S}}_{\epsilon-C_0e^{-R}}^v)$. Thus
        \begin{equation}
            \Vol(G^u\cap(\hat{\mathcal{S}}_{\epsilon+C_1\zeta}^v-\hat{\mathcal{S}}_{\epsilon-C_0e^{-R}}^v))=V(W_v(u,m_v(u)),\epsilon+C_1\zeta)-V(W_v(u,m_v(u)),\epsilon-C_0e^{-R}).
        \end{equation}
        The right hand side of the above equation reaches its maximum when $f_v(u,m_v(u))=0$, i.e., the distance between two centers is $0$; then the maximum is the difference of the volume of two concentric balls. Hence
        when $R$ is sufficiently large, we have
        \begin{align}
            \Vol(G^u\cap(\hat{\mathcal{S}}_{\epsilon+C_1\zeta}^v-\hat{\mathcal{S}}_{\epsilon-C_0e^{-R}}^v))&\ltfrown_n (\epsilon+C_1\zeta)^{k}-(\epsilon-C_0e^{-R})^{k} \notag\\
            &\ltfrown_n \epsilon^{k-1} \zeta.\label{eq:case2}
         \end{align}
    \end{itemize}

    Hence by Proposition \ref{prop:feet_distance}, and \eqref{eq:case2}, when $R$ is sufficiently large, we have
    \begin{equation}
        \begin{aligned}
            \Vol(\hat{\mathcal{S}}_{\epsilon+C_1\zeta}^v-\hat{\mathcal{S}}_{\epsilon-C_0e^{-R}}^v)& \ltfrown_n  (\Vol(b(\hat{\mathcal{S}}_{\epsilon+C_1\zeta}^v))-\Vol(b(\hat{\mathcal{S}}_{\epsilon-C_0e^{-R}}^v)))\cdot\epsilon^{\frac{k-1}{2}}\zeta^{\frac{k+1}{2}}\\
            &\ \ + \Vol(b(\hat{\mathcal{S}}_{\epsilon-C_0e^{-R}}^v))\cdot\epsilon^{k-1} \zeta\\
            &\ltfrown_n \Vol(S^{n-2})\cdot \epsilon^{\frac{k-1}{2}}\zeta^{\frac{k+1}{2}} + \epsilon^{n-2}\cdot \epsilon^{k-1} \zeta\\
            &\asymp \epsilon^{(n^2-n-4)/2}\zeta.
        \end{aligned}
    \end{equation}
    Then by Proposition \ref{mu estimate},
    there exists a constant $M_1$ which only depends on $n$, such that 
    \begin{equation}\label{eq:S-estimate}
        \Vol(\hat{\mathcal{S}}_{\epsilon+C_1\zeta}^v)
        \leq(1+M_1\zeta/\epsilon)
        \Vol(\hat{\mathcal{S}}_{\epsilon-C_0e^{-R}}^v).
    \end{equation}

    By \eqref{eq:leftside}, \eqref{eq:rightside}, \eqref{eq:L-estimate}, and \eqref{eq:S-estimate}, let $L=2M_1+2$, and then \eqref{eq:fiberwise-inequality} is proved.
\end{proof}

\section{Matching and the main theorem} \label{sec:final-proof}

In this last section, we will match each pair of pants to well-matched pants along each boundary component to obtain an incompressible surface in $M$.

\subsection{Matching pattern and the incompressibility} \label{sec:matching_pattern}

We first provide some preliminaries on matching that we will use to complete the proof of the main theorem. 

We start to introduce how we match pants together through average feet, which is a little bit different from the pattern used in \cite{KM12b} and \cite{Ham15} where they used the feet of short orthogeodesics between cuffs. We then will show that if a closed surface in $M$ is made by good pants in this pattern, it is incompressible by applying the definition and the theorem in \cite{Ham15}.

For an oriented closed geodesic $\gamma$ in $M$, we introduce two actions on $N^1(\gamma)$. Let ${\parallel}_\gamma$ be the parallel transport along the orientation of $\gamma$ by distance $1$ and $A_\gamma$ be given by the antipodal map on the unit normal sphere at each point of $\gamma$. We notice that ${\parallel}_\gamma$ commutes with $A_\gamma$, and let 
\begin{equation}
    \tau_\gamma={\parallel}_\gamma \circ A_\gamma=A_\gamma \circ {\parallel}_\gamma\from N^1(\gamma)\to N^1(\gamma).
\end{equation}

Suppose that $P_1,P_2$ are two oriented $(R,\epsilon)$-cuff-good and good pants in $M$ that share a common cuff $\gamma$, for $\epsilon>0$ sufficiently small and $R$ sufficiently large. We then can define:

\begin{defn}\label{matching pattern}
    For $\sigma>0$, $P_1$ and $P_2$ are $\sigma$-\emph{well-matched}, if the following are satisfied:
    \begin{enumerate}
        \item The induced orientations of $\gamma$ as a cuff of $P_1$ and $P_2$ are opposite.

        \item Let $a_1$ be an average foot of $P_1$, then there is an average foot $a_2$ of $P_2$ such that the distance between $\tau_\gamma(a_1)$ and $a_2$ in $N^1(\gamma)$ is less than $\sigma$.
    \end{enumerate}
\end{defn}

\begin{rem}
    We need to clarify that when considering the action of ${\parallel}_\gamma$ on $a_1$, the orientation of $\gamma$ is induced by the orientation of $P_1$. Then this definition is symmetric in $P_1$ and $P_2$, and $\gamma$ is not necessarily oriented.
\end{rem}

We then recall Hamenst{\"a}dt's definition  and result.

\begin{defn}[\cite{Ham15}, Definition 6.1]
    For $\sigma>0$, two pants $P_1,P_2$ are $\sigma$-well attached along a common boundary geodesic $\gamma$ if the following holds true:
    \begin{enumerate}
        \item The orientation of $\gamma$ as a boundary geodesic of $P_1$ and $P_2$ are opposite.

        \item Let $x$ be an endpoint of a short orthogeodesic of $P_1$ on $\gamma$. Then there is an endpoint $y$ of a short orthogeodesic of $P_2$ on $\gamma$ whose oriented distance to $x$ is contained in the interval $[1-\sigma,1+\sigma]$.

        \item Let $v_1$ be the direction of the short orthogeodesic of $P_1$ at $x$ and let $v_2$ be the direction of a short orthogeodesic of $P_2$ at $y$; then the angle between $v_2$ and the image of $-v_1$ under parallel transport along the oriented subarc of $\gamma$ connecting $x$ to $y$ is at most $\sigma$.
    \end{enumerate}
\end{defn}

\begin{prop}[\cite{Ham15}, Proposition 6.2]\label{Ursula matching}
    For any $b>1$, there are $\delta\in(0,\pi/4]$ and $R_0>10$ with the following property. Let $R>R_0$ and let $S\subset M$ be a piecewise immersed closed surface composed of finitely many $(R,\delta)$-good pants which are $R^{-b}$-well attached along their common boundary geodesics. Then $S$ is incompressible.
\end{prop}

\begin{rem}
    In \cite{Ham15}, for a closed oriented surface $S$ of genus $g\ge2$ and a closed Riemannian manifold $M$ of non-positive sectional curvature, a continuous map $f\colon S\to M$ is called \emph{incompressible} if $f_*\colon \pi_1(S)\to\pi_1(M)$ is injective.
\end{rem}

Now we can state our matching result and prove it.

\begin{prop}\label{matching result}
    There are $\delta\in(0,\pi/4]$ and $R_0>10$ such that when $R>R_0$, the following holds: Suppose $S\subset M$ is a piecewise immersed closed surface composed of finitely many $(R,\delta)$-good pants which are $R^{-2}$-well-matched along their common boundary geodesics. Then $S$ is a $\pi_1$-injective surface.
\end{prop}

\begin{proof}
    We choose $b=3/2$ in Proposition \ref{Ursula matching}, then it is sufficient to show that when $R_0$ is large enough, if two $(R,\delta)$-good pants are $R^{-2}$-well-matched, then they are $R^{-3/2}$-well attached. And this is a quick result from the estimate of the difference between average feet and short feet in Proposition \ref{ave close to short}.
\end{proof}

\subsection{Matching pants} \label{sec:matching_pants}

We follow the idea in \cite{KW21} to use Hall's Marriage Theorem to match pants along each good curve. To be more specific, fix an arbitrary good curve $\gamma_0$, and for each oriented good pants $\Pi$ with $\gamma_0$ as a cuff, we want to match $\Pi$ with a well-matched good pants that induces the opposite orientation on $\gamma_0$. We then form a bipartite graph, whose vertices are oriented good pants $P$ with an assigned average foot $a$. There is an edge that joins two vertices if and only if they are well-matched by $\gamma_0$. 

We let $\hat{\Pi}_{R,\epsilon}$ be the set of pairs $(P,a)$ where $P$ is an unoriented good and $(R,\epsilon)$-cuff-good pants in $M$ and $a$ is an average foot of $P$, and $\hat{\Pi}^*_{R,\epsilon}$ is the set of all oriented $(P,a)$.
Since each good pants has three good cuffs and there are two average feet on each good cuff, so the projection $(P,a)\mapsto P$ is a 6-1 mapping. Moreover, $\hat{\Pi}_{R,\epsilon}$ and $\hat{\Pi}^*_{R,\epsilon}$ are finite sets, and $\hat{\Pi}^*_{R,\epsilon}$ is a double cover of $\hat{\Pi}_{R,\epsilon}$ by discarding the orientation. 

For any $(R,\epsilon)$-good curve $\gamma_0$, let $\hat{\Pi}_{R,\epsilon}(\gamma_0)\subset \hat{\Pi}_{R,\epsilon}$ consist of $(P,a)$ where $\gamma_0$ is a cuff of $P$, and $\hat{\Pi}^*_{R,\epsilon}(\gamma_0)\subset \hat{\Pi}^*_{R,\epsilon}$ be the set of oriented $(P,a)$ where $\gamma_0$ is a cuff of $P$.
Let $p:\hat{\Pi}_{R,\epsilon}(\gamma_0)\to N^1(\gamma_0)$ by $p(P,a)=a$, then the point measure induced by $p$ on $N^1(\gamma_0)$ is the actual average foot measure $\nu_a$.

In this subsection, we will prove the following theorem, which is analogous to Theorem 3.1 in \cite{KM12b}, Theorem 5.1 in \cite{KW21}, and Theorem 6.6 in \cite{Rao25}.

\begin{thm}\label{thm:matching}
    For all $\epsilon>0$, there exists $R_0>0$ such that when $R>R_0$, the following holds: Suppose $\gamma_0$ is an $(R,\epsilon)$-good curve in $M$. Then there exists a permutation $\sigma_{\gamma_0}\from \hat{\Pi}_{R,\epsilon}(\gamma_0) \to \hat{\Pi}_{R,\epsilon}(\gamma_0)$, such that for any $\pi\in\hat{\Pi}_{R,\epsilon}(\gamma_0)$,
    \begin{equation}
        d(p(\sigma_{\gamma_0}(\pi)),\tau_{\gamma_0}(p(\pi)))<R^{-2}.
    \end{equation}
\end{thm}

We recall that the Cheeger constant $h(M)$ for a smooth closed $m$-dimensional Riemannian manifold $M$ is defined as
\begin{equation}
    h(M)=\inf\frac{|\partial A|}{\min\{|A|,|M-A|\}},
\end{equation}
among all $m$-dimensional submanifold $A\subset M$. Then for all these $M$, we have the following result.

\begin{thm}[Theorem 5.3 in \cite{KW21}]\label{thm:Cheeger_volume}
    Suppose $A\subset M$ and $|\N_\eta(A)|\leq|M|/2$ for some $\eta>0$, then
    \begin{equation}
        |\N_\eta(A)|\geq (1+\eta h(M))|A|.
    \end{equation}
\end{thm}

To match good pants in the unit normal bundle $N^1(\gamma_0)$, we need to have a lower bound of its Cheeger constant. We can actually have the following general result for an $S^{n-2}$-bundle over a circle.

\begin{prop}\label{prop:Cheeger}
    There exist constants $R_0>0$ such that when $R>R_0$: Suppose 
    $M$ is a $S^{n-2}$-fiber bundle over a round circle with a circumference of $2R$, where $S^{n-2}$ is the standard unit sphere.
    Then $h(M)\geq 1/4 R$.
\end{prop}

\begin{proof}
    Let $S^1_R$ be the circle with a circumference of $2R$. Without loss of generality, we rescale the size of $S^{n-2}$ so that its volume equals 1. Suppose $A$ is a smooth $(n-1)$-dimensional submanifold of $M$ with $|A|\leq |M|/2$. For any $x\in S^1_R$, let $S^{n-2}_x$ be the fiber over $x$ in $M$, and $q_x(A)=|A\cap S^{n-2}_x|$ be the $(n-2)$-dimensional volume of the intersection. Let $h$ be the Cheeger constant of $S^{n-2}$, Then we have
    \begin{equation}
        |\partial A|\geq \int_{x\in S^1_R} h\cdot \min\{q_x(A),|S^{n-2}|-q_x(A)\}\,dx.
    \end{equation}
    Since $M$ is a fiber bundle over a circle, so we also have
    \begin{equation}
        |\partial A|\geq \sup q_x(A) - \inf q_x(A).
    \end{equation}

    If $\sup q_x(A) - \inf q_x(A)\geq 1/4$, then
    \begin{equation}\label{eq:Cheeger_bundle}
        \frac{|\partial A|}{|A|}\geq \frac{1/4}{|A|} \geq \frac{1/2}{|M|} = \frac{1}{4 R}.
    \end{equation}

    If $\sup q_x(A) - \inf q_x(A)\leq 1/4$, then by $|A|\leq |M|/2$, we have $\inf q_x(A)\leq |S^{n-2}|/2=1/2$. Thus $q_x(A)\leq \sup q_x(A)\leq 3/4$, and then $|S^{n-2}|-q_x(A) \geq 1/4$. Therefore
    \begin{equation}
        \begin{split}
            |\partial A|&\geq \int_{x\in S^1_R} h\cdot \min\{q_x(A),|S^{n-2}|-q_x(A)\}\,dx\\
            &\geq \int_{x\in S^1_R} h\cdot \min\{q_x(A),1/4\}\,dx\\
            &\geq \int_{x\in S^1_R} h\cdot \min\{q_x(A),q_x(A)/3\}\,dx\\
            &=\frac{h}{3} \int_{x\in S^1_R} q_x(A)\, dx = \frac{h|A|}{3}.
        \end{split}
    \end{equation}
    Hence $|\partial A|/|A|\geq h/3$.

    Together with \eqref{eq:Cheeger_bundle}, when $R$ is sufficiently large, we have
    \begin{equation*}
        h(M)\geq\min\left\{\frac{1}{4R},\frac{h}{3}\right\}=\frac{1}{4R}.\qedhere
    \end{equation*}
\end{proof}

For any $(R,\epsilon)$-good curve $\gamma_0$, $N^1(\gamma_0)$ is a mapping torus with length $l(\gamma_0)\leq 3R$. Then by Theorem \ref{thm:Cheeger_volume}, we have the following corollary.

\begin{cor}\label{cor:Cheeger}
    If $A\subset N^1(\gamma_0)$ and $|\N_\eta(A)|\leq|N^1(\gamma_0)|/2$ for some $\eta>0$, then
    \begin{equation}
        |\N_\eta(A)|\geq (1+\frac{\eta}{6R})|A|.
    \end{equation}
\end{cor}

Moreover, since $\mu_a$ is quasi-uniform with respect to the Lebesgue measure, we have a similar result for $\mu_a$.

\begin{prop}\label{cor:mu_Cheeger}
    If $A\subset N^1(\gamma_0)$ and $|\N_\eta(A)|\leq|N^1(\gamma_0)|/2$ for some $\eta>0$, then
    \begin{equation}
        \mu_a(\N_\eta(A))\geq (1+\frac{\eta}{6 B^2 R})\mu_a(A),
    \end{equation}
    where $B$ is the constant $B_0$ in Proposition \ref{prop:mu_Lebesgue}.
\end{prop}

\begin{proof}
    By Proposition \ref{prop:mu_Lebesgue} and Corollary \ref{cor:Cheeger}, we have:     
    \begin{equation}         
        \begin{split}             
            \mu_a(\N_\eta(A)-A)
            &\geq B^{-1}
           e^{2R}\epsilon^{(n^2-n+2)/2)}
            |\N_\eta(A)-A|\\
            &\geq B^{-1}e^{2R}\epsilon^{(n^2-n+2)/2)} \frac{\eta}{6R}|A|\\
            &\geq \frac{\eta}{6 B^2 R}\mu_a(A).\qedhere
        \end{split}
    \end{equation}
\end{proof}

Now we can present the proof of Theorem \ref{thm:matching}.

\begin{proof}[Proof of Theorem \ref{thm:matching}]
    Let $\xi=R^{-2}$. By definition, we know that for any $A\subset N^1(\gamma_0)$,
    \begin{equation}
        \nu_a(A)=\#\{\pi\in\hat{\Pi}_{R,\epsilon}(\gamma_0):p(\pi)\in A\}.
    \end{equation}
    Then by Hall's Marriage Theorem, the result follows from the claim that
    \begin{equation}
        \nu_a(\N_{\xi}(A))\geq \nu_a(\tau_{\gamma_0}(A)),
    \end{equation}
    for any $A\subset N^1(\gamma_0)$. Now we prove this claim.

    If $|\N_{\xi/2}(A)|\leq|N^1(\gamma_0)|/2$, by Proposition \ref{maincounting}, Proposition \ref{cor:mu_Cheeger}, and the invariance of $\mu_a$ under $\tau_{\gamma_0}$, we have when $R$ is sufficiently large,
    \begin{equation}
        \begin{aligned}
            \nu_a(\N_{\xi}(A))&\geq \frac{1}{C\Vol(M)}(1-L\zeta/\epsilon)\mu_a(\N_{\xi-\zeta}(A))\\
            &\geq \frac{1}{C\Vol(M)}(1-L\zeta/\epsilon)\mu_a(\N_{\xi/2}(A))\\
            &\geq \frac{1}{C\Vol(M)}(1-L\zeta/\epsilon)(1+\frac{\xi/2-\zeta}{6 B^2 R})\mu_a(\N_{\zeta}(A))\\
            &\geq \frac{1}{C\Vol(M)}(1+L\zeta/\epsilon)\mu_a(\N_\zeta(A))\\
            &= \frac{1}{C\Vol(M)}(1+L\zeta/\epsilon)\mu_a(\N_\zeta(\tau_{\gamma_0}(A)))\\
            &\geq \nu(\tau_{\gamma_0}(A)).
        \end{aligned}
    \end{equation}

    Otherwise, let $A'=N^1(\gamma_0)-\N_\xi(A)$, then $|\N_{\xi/2}(A')|\leq|N^1(\gamma_0)|/2$. Then by the same reasoning, we have $\nu_a(\N_{\xi}(\tau_{\gamma_0}(A')))\geq \nu_a(A')$. Therefore
    \begin{equation}
        \begin{split}
            \nu_a(N^1(\gamma_0))-\nu_a(\tau_{\gamma_0}(A))&=
            \nu_a(\tau_{\gamma_0}(N^1(\gamma_0)-A))\\
            &\geq \nu_a(\tau_{\gamma_0}(\N_\xi(A')))\\
            &\geq \nu_a(N^1(\gamma_0)-\N_\xi(A)).
        \end{split}
    \end{equation}
    Hence $\nu_a(\N_{\xi}(A))\geq \nu_a(\tau_{\gamma_0}(A))$.   
\end{proof}

\subsection{Proof of the main theorem} \label{sec:proof-main}

\begin{proof}[Proof of Theorem \ref{main thm}]
    We want to prove the theorem by showing that, for $n \ge 4$, any closed hyperbolic $n$-manifold $M$ has an immersed incompressible closed surface.

    For any given $\epsilon>0$, we choose $R$ to be sufficiently large such that all previous statements hold for $R$ and $\epsilon$. 

    For any $(R,\epsilon)$-good curve $\gamma_0$, let $\tau_{\gamma_0}$ be a permutation of $\hat{\Pi}_{R,\epsilon}(\gamma_0)$ in Theorem \ref{thm:matching}. Then we apply the doubling trick on $\hat{\Pi}_{R,\epsilon}(\gamma_0)$ and $\tau_{\gamma_0}$. We notice that $2\hat{\Pi}_{R,\epsilon}(\gamma_0)=\hat{\Pi}^*_{R,\epsilon}(\gamma_0)$. Therefore we find an involution $\tau^*$ of $\hat{\Pi}^*_{R,\epsilon}(\gamma_0)$ and match each pants with its image under $\tau^*$. Then we know that $P$ and $P'$ are glued along $\gamma_0$ only if they are $R^{-2}$-well-matched.
    

    Since the number of $(R,\epsilon)$-good curve in $M$ is finite (by Margulis argument), we can repeat the above step for each $(R,\epsilon)$-good curve. Finally, we glue all pants in $\hat{\Pi}^*_{R,\epsilon}$ together and get a closed surface. We notice that this surface may have many connected components, but we can just choose an arbitrary one get a connected closed surface. Then by Proposition \ref{matching result}, this surface will be incompressible in $M$.

\end{proof}

\bibliographystyle{alpha}
\bibliography{reference}

@article{KW21,
  title={\emph{Nearly {Fuchsian} surface subgroups of finite covolume {Kleinian} groups}},
  author={Kahn, Jeremy and Wright, Alex},
  journal={Duke Mathematical Journal},
  volume={170},
  number={3},
  pages={503--573},
  year={2021},
  publisher={Duke University Press}
}

@article{KM12b,
  title={\emph{Immersing almost geodesic surfaces in a closed hyperbolic three manifold}},
  author={Kahn, Jeremy and Markovi\'c, Vladimir},
  journal={Annals of Mathematics},
  pages={1127--1190},
  year={2012},
  publisher={JSTOR}
}

@article{KM15,
title={\emph{The good pants homology and the {Ehrenpreis} conjecture}},
  author={Kahn, Jeremy and Markovi\'c, Vladimir},
  journal={Annals of Mathematics},
  pages={1--72},
  year={2015},
  publisher={JSTOR}
}

@article{CF19,
  title={\emph{Ubiquitous quasi-{Fuchsian} surfaces in cusped hyperbolic 3--manifolds}},
  author={Cooper, Daryl and Futer, David},
  journal={Geometry \& Topology},
  volume={23},
  number={1},
  pages={241--298},
  year={2019},
  publisher={Mathematical Sciences Publishers}
}

@article{LM15,
  title={\emph{Homology of curves and surfaces in closed hyperbolic 3-manifolds}},
  author={Liu, Yi and Markovi\'c, Vladimir},
  journal={Duke Mathematical Journal},
  volume={164},
  number={14},
  pages={2723--2808},
  year={2015},
  publisher={Duke University Press}
}

@article{MZ08,
  title={\emph{Closed quasi-Fuchsian surfaces in hyperbolic knot complements}},
  author={Masters, Joseph D and Zhang, Xingru},
  journal={Geometry \& Topology},
  volume={12},
  number={4},
  pages={2095--2171},
  year={2008},
  publisher={Mathematical Sciences Publishers}
}

@article{BC15,
  title={\emph{Finite-volume hyperbolic 3-manifolds contain immersed quasi-Fuchsian surfaces}},
  author={Baker, Mark D and Cooper, Daryl},
  journal={Algebraic \& Geometric Topology},
  volume={15},
  number={2},
  pages={1199--1228},
  year={2015},
  publisher={Mathematical Sciences Publishers}
}

@article{MZ09,
  title={\emph{Quasi-Fuchsian surfaces in hyperbolic link complements}},
  author={Masters, Joseph D and Zhang, Xingru},
  journal={arXiv preprint arXiv:0909.4501},
  year={2009}
}

@article{Ham15,
  title={\emph{Incompressible Surfaces in Rank One Locally Symmetric Spaces}},
  author={Hamenst{\"a}dt, Ursula},
  journal={Geometric and Functional Analysis},
  volume={25},
  number={3},
  pages={815--859},
  year={2015},
  publisher={Springer}
}

@article{KLM24,
  title={\emph{Surface groups in uniform lattices of some semi-simple groups}},
  author={Kahn, Jeremy and Labourie, Fran{\c{c}}ois and Mozes, Shahar},
  journal={Acta Mathematica},
  volume={232},
  number={1},
  pages={79--220},
  year={2024},
  publisher={International Press of Boston}
}

@misc{KW19,
  title={\emph{Counting connections in a locally symmetric space}},
  author={Kahn, Jeremy and Wright, Alex},
  howpublished = {\url{https://www.math.brown.edu/jk17/Connections.pdf}},
  year={2019}
}

@article{Rao25,
  title={\emph{Subgroups of genus-2 quasi-Fuchsian groups and cocompact Kleinian groups}},
  author={Rao, Zhenghao},
  journal={Geometry \& Topology},
  volume={29},
  number={1},
  pages={495--548},
  year={2025},
  publisher={Mathematical Sciences Publishers}
}

@article{TX18,
  title={\emph{Carrier graphs for representations of the rank two free group into isometries of hyperbolic three space}},
  author={Tan, Ser Peow and Xu, Binbin},
  journal={arXiv preprint arXiv:1807.07277},
  year={2018}
}

@article{Lac10,
  title={Surface subgroups of Kleinian groups with torsion},
  author={Lackenby, Marc},
  journal={Inventiones mathematicae},
  volume={179},
  number={1},
  pages={175--190},
  year={2010},
  publisher={Springer}
}

\end{document}